\documentclass[12pt]{article}
\usepackage{graphicx,psfrag,epsfig}
\usepackage{amssymb,amsmath,amscd,amsthm}
\usepackage{graphicx,psfrag,epsfig}

\usepackage{graphicx}
\usepackage[active]{srcltx}

\newtheorem{theorem}{Theorem}[section]

\newtheorem{lemma}[theorem]{Lemma}

\date{}

\begin{document}

\date{}
\title{Metastability for Non-Linear Random Perturbations of Dynamical Systems}
\author{ M. Freidlin\footnote{Dept of Mathematics, University of Maryland,
College Park, MD 20742, mif@math.umd.edu}, L.
Koralov\footnote{Dept of Mathematics, University of Maryland,
College Park, MD 20742, koralov@math.umd.edu}
} \maketitle

\begin{abstract}
In this paper we describe the long time behavior of  solutions to
quasi-linear parabolic equations with a small parameter at the
second order term and the long time behavior of the corresponding
diffusion processes.
\end{abstract}

{2000 Mathematics Subject Classification Numbers: 60F10, 35K55.}

{ Keywords: Nonlinear Perturbations,  Metastability.}

\section{Introduction}

Consider a dynamical system
\begin{equation} \label{dsyst}
\dot{X}^{x}_t = b( X^{x}_t),~~X^x_0 = x \in \mathbb{R}^d,
\end{equation}
together with its stochastic perturbations
\begin{equation} \label{perturb}
d X^{x,\varepsilon}_t = b(  X^{x,\varepsilon}_t) d t+\varepsilon
\sigma ( X^{x,\varepsilon}_t) d W_t,~~X^{x,\varepsilon}_0 = x \in
\mathbb{R}^d.
\end{equation}
Here $\varepsilon > 0$ is a small parameter, $W_t$ is a Wiener
process in $ \mathbb{R}^d$, and the coefficients $\sigma$ and~$b$
are assumed to be Lipschitz continuous. The diffusion matrix $a(x)
= (a_{ij}(x)) = \sigma(x) \sigma^*(x)$ is assumed to be uniformly
positive definite.

Together with (\ref{perturb}), we can consider the corresponding
Cauchy problem
\begin{equation} \label{eq1cc}
\frac{\partial u^\varepsilon(t,x)}{\partial t}  = L^\varepsilon
u^\varepsilon := \frac{\varepsilon^2}{2} \sum_{i,j=1}^d a_{ij}(x)
\frac{ \partial^2 u^\varepsilon (t,x)}{\partial x_i \partial x_j}
+ b(x) \cdot \nabla_x u^\varepsilon (t,x),~~x \in \mathbb{R}^d,~~t
> 0,
\end{equation}
\begin{equation} \label{bccc}
u^\varepsilon(0,x) = g(x),~x \in \mathbb{R}^d,
\end{equation}
where $g$ is a bounded continuous function.

Suppose for a moment that the vector field $b$ has just one
asymptotically stable equilibrium point $O$  such that all the
points get attracted to $O$ and  $(b(x), x - O) \leq -c|x-O|$ for
some positive constant $c$ and all sufficiently large $|x|$. Then
it is easy to check that
\[
\lim_{(\varepsilon,t) \rightarrow (0,\infty)}  \mathrm{P}
(X^{x,\varepsilon}_t \in U)  = 1
\]
for any neighborhood $U$ of the equilibrium $O$. Taking into
account that the solution $u^\varepsilon$ of
(\ref{eq1cc})-(\ref{bccc}) can be written in the form
$u^\varepsilon(t,x)  = \mathrm{E}g(X^{x,\varepsilon}_t)$ and the
continuity of $g$, we conclude that
\[
\lim_{(\varepsilon,t) \rightarrow (0,\infty)}  u^\varepsilon(t,x)
= g(O).
\]
A similar result holds in the case of a unique compact global
attractor if the system (\ref{dsyst}) has a unique normalized
invariant measure on the attractor. This is the case, for example,
if the system (\ref{dsyst}) in $ \mathbb{R}^2$ has a unique limit
cycle attracting all the trajectories except the unstable
equilibrium inside the cycle.

The situation becomes more complicated if the dynamical system has
more than one asymptotically stable attractor. Assume, for
brevity, that all the attractors are equilibriums $O_1,...,O_n$.
Let $D_i$ be the basin of $O_i$, $1 \leq i \leq n$, and assume
that the set $ \mathbb{R}^d \setminus (D_1 \cup ... \cup D_n)$
belongs to a finite union of surfaces of dimension $d-1$. The long
time behavior of $X^{x,\varepsilon}_t$ and $u^\varepsilon(t,x)$ is
now determined by the transitions of $X^{x,\varepsilon}_t$ between the attractors
$O_1,...,O_n$. These transitions are described by the large
deviation theory for stochastic perturbations of dynamical systems
developed in the late 1960-s (see \cite{FW} and references there).
In particular, the weak limit $\mu$ of the invariant measure
$\mu^\varepsilon$ of the family of processes (\ref{perturb}) was
found. In the generic case, the limiting measure $\mu$ is
concentrated on one of the attractors, which will be denoted by
$O^*$. Then
\[
\lim_{\varepsilon \downarrow 0} \lim_{t \rightarrow \infty}
u^\varepsilon(t,x) = g(O^*).
\]
However, in the case of many attractors, the limiting behavior of
$X^{x,\varepsilon}_t$ and $u^\varepsilon(t,x)$ as $\varepsilon
\downarrow 0$ and $t \rightarrow \infty$ depends on the way in
which $(\varepsilon, t)$ approaches $(0, \infty)$. Roughly
speaking, under natural additional assumptions, there exist a
finite number of regions in the neighborhood of $(0, \infty)$ such
that the limiting distribution of $X^{x,\varepsilon}_t$ and the
limit of  $u^\varepsilon(t,x)$ exist if $(\varepsilon, t)$
approaches $(0, \infty)$ while staying inside one region. For
different regions, these limits are, in general, different.

The corresponding theory of metastability (of sublimiting
distributions) was developed in \cite{F77} (see also \cite{FPhD},
\cite{FW}, \cite{OV}). The notion of a hierarchy of cycles, which
is discussed below, was introduced there. Let $S_{0, T}(\varphi)$
be the action functional for the family  $ X^{x,\varepsilon}_t$ in
$C([0,T], \mathbb{R}^d)$ as $\varepsilon \downarrow 0$
 (\cite{FW}):
\[
S_{0, T}(\varphi) = \frac{1}{2}\int_{0}^{T} \sum_{i,j = 1}^d
a^{ij}(\varphi_t)(\dot{\varphi}^i_t -
b_i(\varphi_t))(\dot{\varphi}^j_t - b_j(\varphi_t)) d t,~~T \geq
0,~ \varphi \in C([0,T], \mathbb{R}^d)
\]
for absolutely continuous $\varphi$, $ S_{0, T}(\varphi) =
+\infty$ for $\varphi$ that are not absolutely continuous. Here
$a^{ij}$ are the elements of the inverse matrix, that is $a^{ij} =
(a^{-1})_{ij}$.  The quasi-potential is defined as
\[
V(x,y) =  \inf_{T, \varphi} \{ S_{0,T}(\varphi): \varphi \in
C([0,T], \mathbb{R}^d), \varphi(0) = x, \varphi(T) = y \},~~x,y
\in \mathbb{R}^d.
\]
Note that while the term ``quasi-potential" is normally applied to
the function $V$ of the variable $y$ with $x$ being a fixed
equilibrium point, we use the same term for the function of two
variables. The hierarchy of cycles is determined by the numbers
\[
V_{ij} = V(O_i, O_j),~~1 \leq i,j \leq n.
\]
The equilibriums $O_1,...,O_n$ are the cycles of rank zero. In the
generic case, for each $O_i$ there exists a unique ``next"
equilibrium $O_l = \mathcal{N}(O_i)$ defined by $V_{i l} = \min_{k:k
\neq i} V_{ik}$. For each sufficiently small $\delta > 0$, with
probability close to one as $\varepsilon \downarrow 0$, the
process $ X^{x,\varepsilon}_t$ that starts in a
$\delta$-neighborhood of $O_i$ will enter a $\delta$-neighborhood
of $\mathcal{N}(O_i)$ before visiting the basins of any of the
equilibriums other than $O_i$ and $\mathcal{N}(O_i)$. The time
before the process enters the neighborhood of $O_l =
\mathcal{N}(O_i)$ is logarithmically equivalent to
$\exp(V_{il}/\varepsilon^2)$. If the sequence $O_i$,
$\mathcal{N}(O_i)$, $\mathcal{N}^2(O_i) = \mathcal{N}(
\mathcal{N}(O_i)) ,...,\mathcal{N}^n(O_i),...$ is periodic, that
is $\mathcal{N}^n(O_i) = O_i$ for some $n$, then a cycle of rank
one appears.  It contains the cycles of rank zero $O_i,
\mathcal{N}(O_i),...,\mathcal{N}^{n-1}(O_i)$. If
$\mathcal{N}^n(O_i) \neq O_i$ for any $n \geq 1$, we say that
$O_i$ forms a cycle of rank one. The entire set of equilibriums is
decomposed into cycles of rank one, which will be denoted by
$C^1_1,...,C^1_{m_1}$. Note that some of the cycles of rank one
may consist of one cycle of rank zero.

Next, the transitions between cycles of rank one can be
considered. Namely, in the generic case, for each cycle $C^1_i$
there is a different cycle $\mathcal{N}(C^1_i)$ of rank one
determined by $V_{ij}$, $1 \leq i,j \leq n$, with the following
property: if the process starts at one of the equilibrium points
in $C^1_i$, then, with probability close to one as $\varepsilon
\downarrow 0$, it will enter a $\delta$-neighborhood of one of the
equilibrium points inside the cycle $\mathcal{N}(C^1_i)$ before
visiting basins of any of the equilibriums outside $C^1_i$ and
$\mathcal{N}(C^1_i)$. This leads to the decomposition of the set
of cycles or rank one into cycles of rank two.

This procedure can be continued inductively until we arrive at a
single cycle of finite rank $R$ which contains all the equilibrium
points. The cycles of rank $r \leq R$ will be denoted by
$C^r_1,...,C^r_{m_r}$.

Let $T^\varepsilon(\lambda) = \exp(\lambda/\varepsilon^2)$. (The
results stated in the paper also hold for $T^\varepsilon(\lambda)
\asymp \exp(\lambda/\varepsilon^2)$, that is if $\varepsilon^2 \ln
T^\varepsilon(\lambda) \rightarrow  \lambda$ as $\varepsilon
\downarrow 0$.) In the generic case, there is a finite set
$\Lambda \subset (0,\infty)$ such that for each $x \in D_1 \cup
... \cup D_n$ and each $\lambda \in (0,\infty) \setminus \Lambda$,
one equilibrium $O_{M(x,\lambda)}$ is defined such that the
measures $\mu^\varepsilon(\Gamma) =
\mathrm{P}(X^{x,\varepsilon}_{T^\varepsilon(\lambda)}\in \Gamma)$
converge weakly to the $\delta$-measure concentrated at
$O_{M(x,\lambda)}$.
The state $O_{M(x,\lambda)}$ is called the metastable state for
the initial point $x$ and the time scale $T^\varepsilon(\lambda)$.

In this paper, instead of the linear problem
(\ref{eq1cc})-(\ref{bccc}), we will consider the Cauchy problem
for the quasi-linear equation with a small parameter
\begin{equation} \label{eq1ccnl}
\frac{\partial u^\varepsilon(t,x)}{\partial t}  = L^\varepsilon
u^\varepsilon := \frac{\varepsilon^2}{2} \sum_{i,j=1}^d
a_{ij}(x,u^\varepsilon) \frac{ \partial^2 u^\varepsilon
(t,x)}{\partial x_i
\partial x_j} + b(x) \cdot \nabla_x u^\varepsilon (t,x),~~x \in
\mathbb{R}^d,~~t
> 0,
\end{equation}
\begin{equation} \label{bcccnl}
u^\varepsilon(0,x) = g(x),~x \in \mathbb{R}^d.
\end{equation}

Equations with diffusion coefficients depending on particle concentration
arise naturally in many applications, in particular in population genetics.
The situation when the drift
$b$ depends on both $x$ and $u^\varepsilon$, with certain
additional assumptions, can also be considered, but we
assume here that $b$ depends only on $x$ for the sake of simplicity.

We assume that the coefficients of equation (\ref{eq1ccnl}) are
Lipschitz continuous and bounded; the matrix $(a_{ij}(x,u))$ is
assumed to be uniformly positive definite. Under these conditions,
problem (\ref{eq1ccnl})-(\ref{bcccnl}) has a unique solution for
any continuous bounded $g(x)$ (see, for instance, \cite{LU}).

A family of processes $X^{x,\varepsilon}_t$, $x \in \mathbb{R}^d$,
satisfying equation (\ref{perturb}) corresponds to each linear
operator $L^\varepsilon$ defined by (\ref{eq1cc}). In the
nonlinear case, a family of processes corresponds to the initial
value problem (\ref{eq1ccnl})-(\ref{bcccnl}). Namely, taking into
account the representation of the solution of the (linear) Cauchy
problem as the expected value of an appropriate functional of the
process, the family corresponding to the problem
(\ref{eq1ccnl})-(\ref{bcccnl}) is defined by the following system
(see \cite{F85}, Ch. 5):
\begin{equation} \label{smtttnew}
d X^{t,x,\varepsilon}_s = b( X^{t,x,\varepsilon}_s) d s +
\varepsilon \sigma(X^{t,x,\varepsilon}_s, u^\varepsilon(t-s,
X^{t,x,\varepsilon}_s)) d W_{s},~~s \leq  t,~~~
X^{t,x,\varepsilon}_0 = x,
\end{equation}
\begin{equation} \label{printe}
u^\varepsilon(t,x) = \mathrm{E}g(X^{t,x,\varepsilon}_{t}),
\end{equation}
where the entries $\sigma_{ij}$ of the matrix $\sigma(x,u)$  are
Lipschitz continuous and $\sigma \sigma^* = a$. The process
$X^{t,x,\varepsilon}_s$ can be viewed as a nonlinear stochastic
perturbation of the dynamical system (\ref{dsyst}).
%

Under the above assumptions on the coefficients and the function
$g$, the solution of the system (\ref{smtttnew})-(\ref{printe})
exists and is unique. The first initial-boundary value problem for
quasi-linear parabolic equation with a small diffusion and the
exit problem for the corresponding processes were studied in
\cite{FK}. The results of the latter paper will be used here.

While the action functional and the quasi-potential were
determined by the time-independent coefficients in the linear
case, now we will consider a family of action functionals and
corresponding quasi-potentials $V_{ij}(c(\lambda))$, $\lambda >
0$. These will be used for times of order $T^\varepsilon(\lambda)
= \exp(\lambda/\varepsilon^2)$. Namely, we will show that the
solution $u^\varepsilon$ of (\ref{eq1ccnl}), in the time scale
$T^\varepsilon(\lambda)$, is very close to a constant $c(\lambda)$
inside $D_i$. We can then define the action functionals and
$V_{ij}(c(\lambda))$ as in the linear case by substituting the
constant $c(\lambda)$ for the second argument in the diffusion
coefficient in the equation.

The main difficulty is that now the action functional and
quasi-potential evolve in time due their dependence on the
(unknown) solution $u^\varepsilon$. Consider, however, a time
interval $[T^\varepsilon(\lambda-\delta),
T^\varepsilon(\lambda)]$, where $\delta$ is small.
%
As will
be seen, $u^\varepsilon$ typically does not change much in time on
this time interval, and the large deviation theory still applies
without drastic modifications, which allows us to express the
limit of $u^\varepsilon(T^\varepsilon(\lambda),x)$, as
$\varepsilon \downarrow 0$, in terms of the limit of
$u^\varepsilon(T^\varepsilon(\lambda-\delta),x)$ and the functions
$V_{ij}(c(\lambda))$. This is the main idea which will allow us to
study the evolution in $\lambda$ of the limit of
$u^\varepsilon(T^\varepsilon(\lambda),x)$. This, in turn, provides
a description of the behavior of $ X^{T^\varepsilon(\lambda)
,x,\varepsilon}_s$ as $\varepsilon \downarrow 0$.

 We will show that if $\lambda$ is sufficiently large, then
the distribution of $ X^{T^\varepsilon(\lambda)
,x,\varepsilon}_{T^\varepsilon(\lambda)}$, even in a generic case,
converges not necessarily to a $\delta$-measure concentrated at an
equilibrium point, but to a distribution on the set of equilibrium
points. Under some natural assumptions this happens, for example,
in the case of two equilibrium points if $V_{12}(c) = V_{21}(c)$
for some value of $c$.
%
%
Therefore, in the case of nonlinear perturbations, the notion of a
metastable state should be replaced by the notion of a metastable
distribution.

Note that metastable distributions (rather than states) arise also
in the case of linear parabolic equations. For example, if the
non-perturbed system, say in $ \mathbb{R}^2$, has two
asymptotically stable limit cycles attracting the entire space,
other than the separatrices, then each of the invariant
distributions on those cycles will be the metastable distribution
for the appropriate initial states and time scales. Metastable
distributions on an asymptotically stable attractor arise in
physical models (see various models and references in \cite{OV}).
However, in the case considered here, the metastable distributions
are supported on several separated asymptotically stable
attractors. Similar metastable distributions arise also when
perturbations of nearly-Hamiltonian systems are considered (see
\cite{AF}, \cite{F}), but because of different reasons.



Since the quasi-potential changes in time,  the relative stability
of attractors also changes in time, possibly leading to changes in
the hierarchy of cycles. We will mostly be concerned with the
situation when the hierarchy of cycles does not change. This is
the case, for example, if there are only two equilibrium points or
if the matrix $a_{ij}(x,u)$ is close enough to a diffusion matrix
independent of $u$.
%
An example with a change in the hierarchy of cycles is considered
in Section~\ref{change}.

If the system has $n$ asymptotically stable equilibrium points (or
more general stable attractors), the number of different (even
generic) cases which should be considered  grows very fast with
$n$: one should consider not just different hierarchies of cycles,
but also different relations between the values of the initial
function $g$ at the equilibriums and various behaviors of
$V_{ij}(c)$ as $c$ changes. Therefore we consider in more detail
the case of two attractors and describe the result in the case of
three attractors. The general result is not presented, but we
believe that the methodology developed in this paper for the case
of small $n$ works in general (generic) case.

In Section~\ref{defs} we introduce some of the definitions and
discuss the notion of the hierarchy of cycles in more detail. We
also state the lemmas that can be used to describe the long-time
behavior of a process whose time-dependent coefficients are close
to functions that do not depend on time. In
Sections~\ref{example1} and \ref{example2} we consider
 a system with two equilibriums and a system with three
equlibriums on the real line in the case when the hierarchy of
cycles is preserved. In Section~\ref{general} we formulate a
general result for the case when the hierarchy of cycles is
preserved. In Section~\ref{change} we study the asymptotics of the
solution to the parabolic equation for a system in which a
bifurcation in  the hierarchy of cycles occurs.

\section{Notations. Diffusion Processes Corresponding to the Nonlinear Problem}
\label{defs}

 Let $\alpha(x)$ be a symmetric $d \times d$ matrix whose elements
$\alpha_{ij}(x)$ are  Lipschitz continuous  with Lipschitz
constant $L$  and satisfy
\begin{equation} \label{elliptic}
k|\xi|^2 \leq \sum_{i,j =1}^d \alpha_{ij}(x) \xi_i \xi_j \leq
K|\xi|^2,~~ x \in \mathbb{R}^d,~~ \xi \in \mathbb{R}^d.
\end{equation}
 Let
$\alpha^{ij}$ be the elements of the inverse matrix, that is
$\alpha^{ij} = (\alpha^{-1})_{ij}$, and $\sigma$ be a square
matrix such that $\alpha = \sigma \sigma^*$. We choose $\sigma$ in
such a way that $\sigma_{ij}$ are also Lipschitz continuous.

We assume that all the attractors of the  bounded Lipschitz
continuous vector field $b$ are equilibriums $O_1,...,O_n$. Assume
that their domains of attraction $D_1,...,D_n$ are such that the
set $ \mathbb{R}^d \setminus (D_1 \cup ... \cup D_n)$ belongs to a
finite union of surfaces of dimension $d-1$. We also assume that
there are $r > 0$ and $c > 0$ such that
\begin{equation} \label{attra}
(b(x), x-O_i) \leq -c|x-O_i|^2
\end{equation}
whenever $x$ is in the $r$-neighborhood of $O_i$, $1 \leq i \leq
n$.

Let $S^\alpha_{0, T}$ be the normalized action functional for the
family of processes $X^{x,\varepsilon}_t$ satisfying
\begin{equation} \label{bsk}
 d X^{x,\varepsilon}_t = b(  X^{x,\varepsilon}_t) d t + \varepsilon
\sigma ( X^{x,\varepsilon}_t) d W_t,~~X^{x, \varepsilon}_0  = x,
\end{equation}
 where $b$ is  a bounded Lipschitz continuous vector field
on $ \mathbb{R}^d$. Thus
\[
S^\alpha_{0, T}(\varphi) = \frac{1}{2}\int_{0}^{T} \sum_{i,j =
1}^d \alpha^{ij}(\varphi_t)(\dot{\varphi}^i_t -
b_i(\varphi_t))(\dot{\varphi}^j_t - b_j(\varphi_t)) d t
\]
for absolutely continuous $\varphi$ defined on $[0, T]$,
$\varphi_0 = x$,  and $ S^\alpha_{0, T}(\varphi) = \infty$ if
$\varphi$ is not absolutely continuous or if $\varphi_0 \neq x$
(see \cite{FW}). Let $V^\alpha(x,y)$ be the quasi-potential for
the family $X^{x,\varepsilon}_t$ in $\mathbb{R}^d$, that is
\begin{equation} \label{vvv}
V^\alpha(x,y) = \inf_{T, \varphi} \{ S^\alpha_{0,T}(\varphi):
\varphi \in C([0,T], \mathbb{R}^d),  \varphi(0) = x, \varphi(T) =
y \},~~x,y \in \mathbb{R}^d.
\end{equation}
Let $V^\alpha_{i j} = V^\alpha(O_i,O_j)$.
For a given function $\alpha$, we define inductively the following
objects (see \cite{F77}, \cite{FW} for a detailed exposition).

(a) The hierarchy of cycles $C^r_1,...,C^r_{m_r}$, $r \leq R$.

(b) The notion of the ``next" equilibrium $\nu(C^r_i)$ and the
``next" cycle $\mathcal{N}(C^r_i)$ of the same rank for a cycle
$C^r_i$ of rank less than $R$.

(c) The transition rates
$V^{\alpha}_{C^r_i, O_j}$, $1 \leq i \leq m_r$, $1 \leq j \leq n$,
$O_j \notin C^r_i$, from a cycle to equilibriums outside this
cycle.

%


For $r = 0$, we define $C^0_i = \{O_i\}$,
$V^{\alpha}_{C^0_i, O_j} = V^{\alpha}_{i j}$.
Assume that the cycles of rank $r$ and the transition rates from
those cycles to equilibrium points have been defined. We
define~$O_j$ to be the next equilibrium after $C^r_i$ if
$\min_{j: O_j \notin C^r_i} V^{\alpha}_{C^r_i, O_j} $ is achieved
at $j$.
\\

\noindent
 {\bf Assumption A.} The minimum
$\min_{j: O_j \notin C^r_i} V^{\alpha}_{C^r_i, O_j}$ is achieved
for a single value of $j$.
\\

We will write $O_j = \nu(C^r_i)$ to express that $O_j$ is the next
equilibrium after $C^r_i$.
We say that the cycle $C^r_l$ of rank $r$ is the next after
$C^r_i$ if $C^r_l$ contains $\nu(C^r_i)$.  We will express this
relation by writing $C^r_l = \mathcal{N}(C^r_i)$. Starting from a
cycle $C^r_i$ of rank $r$, we can form the sequence $C^r_i,
\mathcal{N}(C^r_i), \mathcal{N}^2 (C^r_i),...$ by using the
operation ``next". If this sequence is periodic, that is $C^r_i =
\mathcal{N}^n(C^r_i)$ for some $n$, then the cycles $C^r_i,...,
\mathcal{N}^{n-1} (C^r_i)$ form a
 cycle of rank $r+1$.  If $C^r_i \neq \mathcal{N}^n(C^r_i)$ for
 any $n \geq 1$, then $C^r_i$ is said to form a cycle of rank
 $r+1$.
This way, the collection of all the cycles of rank $r$ is
decomposed in a  union of non-intersecting cycles of rank $r+1$.

If $C^r_{1},...,C^r_{s}$ form a cycle of rank $r+1$, which will be
denoted by $\Gamma$, we define
$V^{\alpha}_{\Gamma, O_j}$ as
\begin{equation} \label{srkp1}
V^{\alpha}_{\Gamma, O_j} =   \max_{1 \leq m \leq s}
V^\alpha_{C^r_{m}, \nu(C^r_{m})} + \min_{1 \leq m \leq s}
(V^{\alpha}_{C^r_{m}, O_j} - V^\alpha_{C^r_{m},
\nu(C^r_{{m}})}),~~O_j \notin \Gamma.
\end{equation}
%
%
 We can continue this procedure until we
arrive at a single cycle of highest rank $R$.

If $\Gamma$ is a cycle, we define $D_\Gamma = \cup_{i: O_i \in
\Gamma} D_i$. As follows from \cite{F77}, \cite{FW}, if the
process (\ref{bsk}) starts in $D_\Gamma$, where $\Gamma$ is a
cycle of rank $r < R$,  then with probability which tends to one
as $\varepsilon \downarrow 0$ it will leave $D_\Gamma$ and enter a
small neighborhood of $\nu(\Gamma)$ in time $T(\varepsilon) \asymp
\exp(V^\alpha_{ \Gamma, \nu(\Gamma)} /\varepsilon^2)$.

 Next we discuss the long-time behavior of processes
whose diffusion coefficients are time-dependent, but are close to
functions that do not depend on time. For $T > 0$ and $\varphi,
\psi \in C([0,T], \mathbb{R}^d)$, we define $\rho_T(\varphi, \psi)
= \sup_{t \in [0,T]}|\varphi(t) - \psi(t)|$.

Let $\widetilde{\alpha}^\varepsilon(t,x)$ be a uniformly positive
definite symmetric $d \times d$ matrix whose elements
$\widetilde{\alpha}_{ij}^\varepsilon$ are continuous in $(t,x)$
and  Lipschitz continuous in $x$.
Let $\widetilde{\sigma}^\varepsilon$ be a square matrix such that
$\widetilde{\alpha}^\varepsilon = \widetilde{\sigma}^\varepsilon
(\widetilde{\sigma}^\varepsilon)^*$. We choose
$\widetilde{\sigma}^\varepsilon$ in such a way that
$\widetilde{\sigma}_{ij}^\varepsilon$ are also  continuous in
$(t,x)$ and  Lipschitz continuous in $x$.

Let $\widetilde{X}^{x,\varepsilon}_t$ satisfy $\widetilde{X}^{x,
\varepsilon}_0 = x$ and
\begin{equation} \label{diffpt}
d \widetilde{X}^{x,\varepsilon}_t = b(
\widetilde{X}^{x,\varepsilon}_t) d t +\varepsilon
\widetilde{\sigma}^\varepsilon (t,
\widetilde{X}^{x,\varepsilon}_t) d W_t,
\end{equation}
 where $b$ is the same as above.
  The law of this process depends on
$\widetilde{\sigma}^\varepsilon$ only through
$\widetilde{\alpha}^\varepsilon = \widetilde{\sigma}^\varepsilon
(\widetilde{\sigma}^\varepsilon)^*$. We will assume that the
diffusion coefficients for the process
$\widetilde{X}^{x,\varepsilon}_t$ are close to those of
${X}^{x,\varepsilon}_t$. Namely, let us assume that
\begin{equation} \label{close}
 \sup_{(t,x) \in \mathbb{R}^+
\times \mathbb{R}^d} | \widetilde{\alpha}_{ij}^\varepsilon(t, x)-
{\alpha}_{ij}(x)| \leq \varkappa,
\end{equation}
where $\varkappa$ is small. The reason to introduce the process
$\widetilde{X}^{x,\varepsilon}_t$ is that we would like to study
the behavior of the process $ X^{t,x,\varepsilon}_s$ given by
(\ref{smtttnew})-(\ref{printe}) on a time interval where the
variable  $u^\varepsilon$ found inside the diffusion coefficient
of (\ref{smtttnew}) does not change much. Since a-priori we don't
know much about the behavior of the diffusion coefficients in
(\ref{smtttnew}) (other than that they don't significantly change
in time on a certain time interval), it is convenient to consider
a generic process whose diffusion coefficients are close to
functions that don't depend on time.

The next two lemmas show that $S^\alpha$ serves a purpose similar
to the action functional for the process
$\widetilde{X}^{x,\varepsilon}_t$, even though the diffusion
coefficients for the process are time-dependent.
\begin{lemma} \label{ler1}
Suppose that $b$ is fixed, $\alpha$ and
$\widetilde{\alpha}^\varepsilon$ are as above, and positive
constants  $k$, $K$ and $L$ are fixed. For any $\delta$, $\gamma$
and $C$ there exist $\varkappa > 0$ and $\varepsilon_0
> 0$ such that
\[
\mathrm{P}(\rho_T(\widetilde{X}^{x,\varepsilon}_t, \varphi) <
\delta) \geq \exp(-\varepsilon^{-2}[S^\alpha_{0,T}(\varphi) +
\gamma])
\]
for $\varepsilon < \varepsilon_0$ and $T > 0$, $\varphi \in
C([0,T], \mathbb{R}^d)$ such that $\varphi(0) = x$ and $T +
S^\alpha_{0,T}(\varphi) < C$.
\end{lemma}
\begin{lemma} \label{ler2} Suppose that $b$ is fixed, $\alpha$ and
$\widetilde{\alpha}^\varepsilon$ are as above, and  the positive
constants  $k$, $K$ and $L$ are fixed. For $x \in \mathbb{R}^d$,
$T >0$ and $s \geq  0$, put
\[
\Phi(s) = \{ \varphi \in C([0,T], \mathbb{R}^d), \varphi(0) = x,
S^\alpha_{0,T}(\varphi) \leq s \}.
\]
For any $T > 0$, $\delta > 0$, $\gamma > 0$ and $s_0 > 0$, there
exist $\varkappa > 0$ and $\varepsilon_0 > 0$ such that for $x \in
\mathbb{R}^d$, $0 < \varepsilon \leq \varepsilon_0$ and $s \leq
s_0$, we have
\[
\mathrm{P}(\rho_T(\widetilde{X}^{x,\varepsilon}_t, \Phi(s)) \geq
\delta) \leq \exp(-\varepsilon^{-2}[s - \gamma]).
\]
\end{lemma}
Note that the choice of $\varkappa$ and $\varepsilon_0$ in the
above lemmas depends on $\alpha$ and
$\widetilde{\alpha}^\varepsilon$ only through $k$, $K$ and $L$.
\\
\\
{\it Sketch of the proof of Lemmas \ref{ler1} and \ref{ler2}.} The
proof of these lemmas is similar to the proof of the fact that
$S^\alpha_{0, T}(\varphi)$ serves as an action functional for the
process $X^{x,\varepsilon}_t$ given in (\ref{bsk}) (see
\cite{FWa}, \cite{DSt}). In order to apply the method based on the
Euler approximations (see Section 1.4 of \cite{DSt}), we need to
show that a process with constant diffusion coefficients is close
to a process with slightly perturbed coefficients in the following
sense:

Let $Y^{x,\varepsilon}_t$, $\widetilde{Y}^{x,\varepsilon}_t$
satisfy
\begin{equation} \label{bskvv}
 d Y^{x,\varepsilon}_t = b d t + \varepsilon
\sigma d W_t,~~Y^{x, \varepsilon}_0  = x,
\end{equation}
\begin{equation} \label{bskww}
 d \widetilde{Y}^{x,\varepsilon}_t = b d t + \varepsilon
(\sigma + \delta(t,\widetilde{Y}^{x,\varepsilon}_t)) d
W_t,~~\widetilde{Y}^{x, \varepsilon}_0  = x,
\end{equation}
where $b$ is a constant vector, $\sigma$ is a constant matrix, and
$\delta(t,x)$ is a matrix whose entries are continuous in $(t,x)$
and Lipschitz continuous in $x$. Then for each positive $h$, $A$
and $T$, there is a positive $\delta_0$ such that
\begin{equation} \label{close3}
\mathrm{P}(\sup_{t \leq T} |Y^{x,\varepsilon}_t -
\widetilde{Y}^{x,\varepsilon}_t| > h) \leq \exp(-A/(\varepsilon^2
T))
\end{equation}
if $\sup_{t,x}||\delta(t,x)|| \leq \delta_0$  and $\varepsilon$ is
sufficiently small. (Here we define $||\delta|| = \sqrt{\sum_{i,j
= 1}^d (\delta^{ij})^2}~$).  To prove (\ref{close3}), we note that
the $i$-th component of the difference satisfies
\[
M^i_t := (Y^{x,\varepsilon}_t - \widetilde{Y}^{x,\varepsilon}_t)^i
= \varepsilon \int_0^t \sum_{j=1}^d
\delta^{ij}(t,\widetilde{Y}^{x,\varepsilon}_s)) d W^j_s.
\]
The right hand side is a martingale with quadratic variation
satisfying
\[
\langle M^i \rangle_t \leq \varepsilon^2
t\sup_{t,x}||\delta(t,x)||^2 \leq \varepsilon^2 t \delta_0^2.
\]
Therefore
\[
\sup_{t \leq T}|M^i_t| \leq \sup_{t  \leq T} |\widetilde{W} (
\varepsilon^2 t  \delta_0^2)|,
\]
where $ \widetilde{W}$ is a standard Brownian motion. Therefore
\[
\mathrm{P}(\sup_{t \leq T} |Y^{x,\varepsilon}_t -
\widetilde{Y}^{x,\varepsilon}_t| > h) \leq d \mathrm{P}(\sup_{t
\leq T} |\widetilde{W} ( \varepsilon^2 t  \delta_0^2)| >
\frac{h}{d}),
\]
which can clearly be made smaller than the right hand side of
(\ref{close3}) by selecting a sufficiently small~$\delta_0$. \qed

We next state a corollary of the above two lemmas that will be
used in the paper. Given a domain $D$ and $\delta >0$, we define
\[D^\delta = \{ x \in D: {\rm dist}(x,
\partial D) \geq \delta,~|x| \leq 1/\delta \}.
\]
Let $x_0$ be an asymptotically stable equilibrium of $b$ and $D$
be a domain attracted to $x_0$.
Let
\[
v = \inf_{T, \varphi} \{ S^\alpha_{0,T}(\varphi): \varphi \in
C([0,T], \overline{D}), \varphi(0) = x_0, \varphi(T) \in \partial
D \}.
\]
\begin{lemma} \label{mlml}
Suppose that $b$ is fixed, $\alpha$ is Lipschitz continuous
 with Lipschitz
constant $L$, $\widetilde{\alpha}^\varepsilon$ is continuous in
$(t,x)$ and  Lipschitz continuous in $x$,  and
\[
k|\xi|^2 \leq \sum_{i,j =1}^d {\alpha}_{ij}( x) \xi_i \xi_j \leq
K|\xi|^2~~ for~~x \in D,~~\xi \in \mathbb{R}^d,
\]
\begin{equation} \label{ellik}
k|\xi|^2 \leq \sum_{i,j =1}^d \widetilde{\alpha}_{ij}^\varepsilon
(t, x) \xi_i \xi_j \leq K|\xi|^2~~ for~~(t,x) \in \mathbb{R}^+
\times D,~~\xi \in \mathbb{R}^d.
\end{equation}
 For each $\delta > 0$ there
are $\varkappa > 0$ and a function $\rho(\varepsilon)$ (that
depend on  $\alpha$ and $\widetilde{\alpha}$ through $L, k$ and
$K$) such that $\lim_{\varepsilon \downarrow 0} \rho(\varepsilon)
= 0$ and
\[
\sup_{ (t,x) \in  [T^\varepsilon(\delta), T^\varepsilon(v-\delta)]
\times D^\varkappa} \mathrm{P}(|\widetilde{X}^{x,\varepsilon}_t -
x_0| < \delta,~\widetilde{X}^{x,\varepsilon}_s \in D~~{\rm for}~~s
\leq t) \geq 1 - \rho(\varepsilon),
\]
provided that
\[
 \sup_{(t,x) \in \mathbb{R}^+
\times D^\varkappa }| \widetilde{\alpha}_{ij}^\varepsilon(t, x)-
{\alpha}_{ij}(x)| \leq \varkappa.
\]
\end{lemma}
This lemma can be easily proved using a modification of Theorems
4.2 and 4.3 from Chapter 4 of \cite{FW} if we substitute our
Lemmas \ref{ler1} and \ref{ler2} for the corresponding results
concerning the case of time-independent coefficients.

The next simple lemma does not require the proximity of
$\widetilde{\alpha}^\varepsilon$ to $\alpha$, but only the
boundedness of the entries of $\widetilde{\alpha}^\varepsilon$. It
can be proved by standard arguments from large deviation theory
(compare with chapter 3 of \cite{FW}).

\begin{lemma} \label{letf} Suppose that $b$ is fixed and $\widetilde{\alpha}^\varepsilon$
is continuous in $(t,x)$ and  Lipschitz continuous in $x$ and
satisfies (\ref{ellik}).
For any compact $M \subset D$, there is $v_0 > 0$ which depends on
$\widetilde{\alpha}^\varepsilon$ only through $K$ such that for
each $\delta \in (0, v_0)$ there is a function $\rho(\varepsilon)$
such that $\lim_{\varepsilon \downarrow 0} \rho(\varepsilon) = 0$
and
\[
\sup_{(t,x) \in  [T^\varepsilon(\delta), T^\varepsilon(v_0)]
\times M} \mathrm{P}(|\widetilde{X}^{x,\varepsilon}_t - x_0| <
\delta,~\widetilde{X}^{x,\varepsilon}_s \in D~~{\rm for}~~s \leq
t) \geq 1 - \rho(\varepsilon).
\]
\end{lemma}

Note that the quasi-potential can be defined by (\ref{vvv}) even
if $\alpha$ has some discontinuities.   We shall be particularly
interested in the structure of the hierarchy of cycles and the
exponential transition times for functions $\alpha$ which are of
the form $\alpha = a(x, f(x))$, where $f$ is constant on each
$D_i$.
The reason for that is that the solution of
(\ref{eq1ccnl})-(\ref{bcccnl}) is nearly constant inside each of the domains
$D_i^\delta = \{ x \in D_i: {\rm dist}(x,
\partial D_i) \geq \delta,~|x| \leq 1/\delta \}$, $\delta > 0$, $1
\leq i \leq n$,  for $\varepsilon$ small enough, as follows from
the following lemma.
\begin{lemma} \label{letwomk} Let  $u^\varepsilon$ be the solution of
 (\ref{eq1ccnl})-(\ref{bcccnl}).
For every positive $\lambda_0$ and $\delta$ there is a positive
$\varepsilon_0$ such that
\begin{equation} \label{kapmk} |u^\varepsilon(T^\varepsilon(\lambda), x) -
u^\varepsilon(T^\varepsilon(\lambda),O_i)| \leq \delta
\end{equation}
whenever $x \in D_i^\delta$, $\varepsilon \leq \varepsilon_0$ and
$\lambda \geq \lambda_0$.
\end{lemma}
For a proof of this lemma we refer the reader to \cite{FK}, where
the same statement was proved in the case of a single domain.

\section{The case of two equilibrium points} \label{example1}

In this section we assume that there are two asymptotically stable
equilibrium points $O_1, O_2 \in \mathbb{R}^d$. Let $D_1 \subset
\mathbb{R}^d$ be the set of points in $\mathbb{R}^d$ which are
attracted to $O_1$ and $D_2 \subset \mathbb{R}^d$ the set of
points attracted to $O_2$. We assume that $D_1 \cup D_2 \in
\mathbb{R}^d \setminus S$, where $S$ is a $(d-1)$-dimensional
manifold. Note that in the case of two equilibrium points,
 the hierarchy of cycles is always the same: $O_1$ and $O_2$ are cycles of rank zero,
and there is one cycle of rank one which contains both $O_1$ and
$O_2$.



Let $g_{\rm min} = \inf_{x \in \mathbb{R}^d} g(x)$ and $g_{\rm
max} = \sup_{x \in \mathbb{R}^d} g(x)$.
Define the functions $M_{12}$, $M_{21}: [g_{\rm min},
g_{\rm max}] \rightarrow \mathbb{R}$ via
\[
M_{12}(c) =  {V}^{ a(\cdot, c)}_{O_1,O_2},~~
M_{21}(c) =  {V}^{ a(\cdot, c)}_{O_2,O_1}.
\]
These functions are shown on Figure 1. It is not difficult to
check that the constant $c$ in the second argument of $a$ can be
replaced by any function equal to $c$ on $D_1$ in the definition
of $M_{12}$ and equal to $c$ on $D_2$ in the definition of
$M_{21}$ without affecting the values of $M_{12}(c)$ and
$M_{21}(c)$.


\begin{figure}[htbp]
 \label{pic2a}
  \begin{center}
    \begin{psfrags}
     \includegraphics[height=6.5in, width= 3.5in,angle=90]{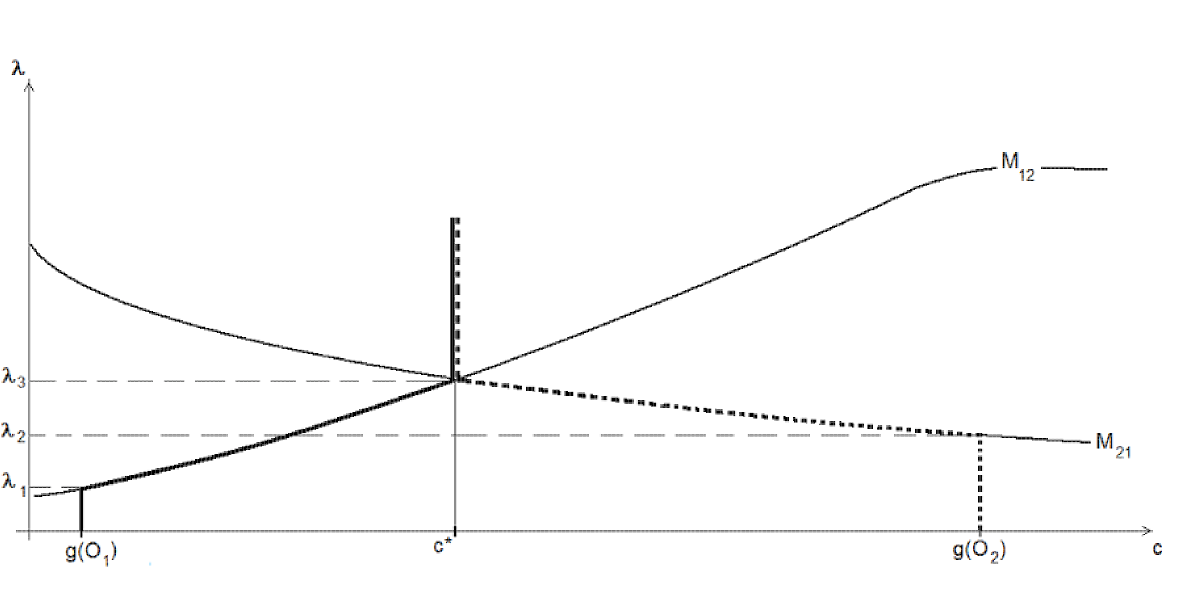}
    \end{psfrags}
    \caption{The case of two equilibrium points}
  \end{center}
\end{figure}

Without loss of generality we may assume that $g(O_1) \leq
g(O_2)$. Let  $\lambda_1=M_{12}(g(O_1))$  and $\lambda_2=
M_{21}(g(O_2))$. In order to formulate the results on the
asymptotics of $u^\varepsilon(T^\varepsilon(\lambda), x)$, we need
the functions $c^1(\lambda)$ and $c^2(\lambda)$, $\lambda > 0$,
defined as follows:

\begin{equation} \label{hhnn1}
c^1(\lambda) =
\left\{\begin{array}{lll}g(O_1),~~~~~~~~~~~~~~~~~~~~~~~~~~~~~~~~~~~~~~~~~~~~~~~~~~~~~~~~~~~~~~~\,0
< \lambda < \lambda_1,\\ \min \{g(O_2), \min\{c: c \in [g(O_1),
g(O_2)], M_{12}(c) = \lambda \} \},~~\,\,\lambda \geq
\lambda_1,\end{array}\right.
\end{equation}
\begin{equation} \label{hhnn2}
c^2(\lambda) =
\left\{\begin{array}{lll}g(O_2),~~~~~~~~~~~~~~~~~~~~~~~~~~~~~~~~~~~~~~~~~~~~~~~~~~~~~~~~~~~~~~~~\,0
< \lambda < \lambda_2,\\ \max\{g(O_1), \max\{c: c \in [g(O_1),
g(O_2)], M_{21}(c) = \lambda \} \},~~\,\,\lambda \geq
\lambda_2.\end{array}\right.
\end{equation}

%
%
%
%

 Let $\lambda_3 = \inf\{\lambda:
c^1(\lambda) \geq c^2(\lambda)\}$. Assume that at least one of the
functions $c^1$ and $c^2$ is continuous at $\lambda_3$. Let $c^* =
c^1(\lambda_3)$ if $c^1$ is continuous at $\lambda_3$ and $c^* =
c^2(\lambda_3)$ otherwise. Let $\overline{c}^1(\lambda) =
\min(c^1(\lambda), c^*)$ and $\overline{c}^2(\lambda) =
\max(c^2(\lambda), c^*)$. On Figure 2, the graphs of
$\overline{c}^1$ and $\overline{c}^2$ are denoted by the thick and
the dotted lines, respectively.

The asymptotics of $u^\varepsilon(T^\varepsilon(\lambda), x)$ is
described by the following theorem. Later, we will use this result
to describe the behavior of the process $X^{t,x,\varepsilon}_s$
when $\varepsilon \downarrow 0$.

\begin{theorem} \label{mtex1}
Let the above assumptions  be satisfied. Suppose that the function
$\overline{c}^1(\lambda)$ is continuous at a point $\lambda \in
(0, \infty)$. Then for every $\delta > 0$ the following limit
\[
\lim_{\varepsilon \downarrow 0}
u^\varepsilon(T^\varepsilon(\lambda), x) = \overline{c}^1(\lambda)
\]
is uniform in $x \in D_1^\delta$. Suppose that the function
$\overline{c}^2(\lambda)$ is continuous at a point $\lambda \in
(0, \infty)$. Then for every $\delta > 0$ the following limit
\[
\lim_{\varepsilon \downarrow 0}
u^\varepsilon(T^\varepsilon(\lambda), x) = \overline{c}^2(\lambda)
\]
is uniform in $x \in D_2^\delta$.
\end{theorem}
\proof  Let us show that if ${c}^1$ is continuous at $\lambda$,
then
\begin{equation} \label{toshow}
\limsup_{\varepsilon \downarrow 0} \sup_{x \in D_1^\delta}
u^\varepsilon(T^\varepsilon(\lambda), x) \leq {c}^1(\lambda).
\end{equation}
Similarly, if  ${c}^2$ is continuous at $\lambda$, then
\begin{equation} \label{toshow2}
\liminf_{\varepsilon \downarrow 0} \inf_{x \in D_2^\delta}
u^\varepsilon(T^\varepsilon(\lambda), x) \geq {c}^2(\lambda).
\end{equation}
Due to Lemma~\ref{letwomk}, in order to prove (\ref{toshow}), it
is sufficient to show that
\begin{equation} \label{twhh}
\limsup_{\varepsilon \downarrow 0}
u^\varepsilon(T^\varepsilon(\lambda), O_1) \leq {c}^1(\lambda),
\end{equation}

Note that by Lemma~\ref{letf} and (\ref{printe}) there is a
positive $v_0$ such that for every $0 < \delta < v_0$ there is
$\varepsilon_0 > 0$ such that
\begin{equation} \label{initxx}
 |u^\varepsilon(T^\varepsilon(\lambda), x) - g(O_i)| \leq \delta
\end{equation}
whenever $x \in D^\delta_i$, $0 < \varepsilon \leq \varepsilon_0$
and $\delta \leq \lambda \leq v_0$.

If (\ref{twhh}) fails for a certain value of $\lambda$, then due
to continuity of the functions $u^\varepsilon(t,O_i)$ in $t$, it
follows from (\ref{initxx}) that for an arbitrarily small $\delta'
> 0$ there are sequences $\varepsilon_n \downarrow 0$ and $\lambda_n \in
[\delta', \lambda]$ such that
\[ u^{\varepsilon_n}(t, O_1) \leq {c}^1(\lambda) +
\delta',~~~ T^{\varepsilon_n}(\delta') \leq t \leq
T^{\varepsilon_n}(\lambda_n)
\]
and
\[
u^{\varepsilon_n}(T^{\varepsilon_n}(\lambda_n), O_1) =
{c}^1(\lambda) + \delta'.
\]
Take $\delta'' \in (0, \delta')$ which will be specified later.
Due to the continuity of $u^{\varepsilon_n}(t,O_1)$ in $t$, we can
find a sequence $\mu_n \in [\delta', \lambda_n)$ such that
\[
u^{\varepsilon_n}(T^{\varepsilon_n}(\mu_n), O_1) = {c}^1(\lambda)
+ \delta''
\]
and
\begin{equation} \label{clll}
u^{\varepsilon_n}(t, O_1) \in [{c}^1(\lambda) + \delta'',
{c}^1(\lambda) + \delta']~~{\rm for}~~t \in [
T^{\varepsilon_n}(\mu_n), T^{\varepsilon_n}(\lambda_n)].
\end{equation}
 We
can express $ u^{\varepsilon_n}(T^{\varepsilon_n}(\lambda_n),
O_1)$ in terms of the process $  X^{
T^{\varepsilon_n}(\lambda_n),O_1,\varepsilon}_s$ and the solution
at the earlier time $ T^{\varepsilon_n}(\mu_n)$ as follows
 \begin{equation}\label{xpr}
u^{\varepsilon_n}(T^{\varepsilon_n}(\lambda_n), O_1) =
\mathrm{E}u^{\varepsilon_n}\left( T^{\varepsilon_n}(\mu_n), X^{
T^{\varepsilon_n}(\lambda_n),O_1,\varepsilon_n}_{T^{\varepsilon_n}(\lambda_n)
- T^{\varepsilon_n}(\mu_n)}\right).
\end{equation}
Since ${c}^1$ is continuous at $\lambda$, there are arbitrarily
small $\delta'>0$ such that $M_{12}({c}^1(\lambda) + \delta')
> M_{12}({c}^1(\lambda)) = \lambda$. Since $\lambda_n \leq \lambda$,
a process starting at $O_1$ and satisfying (\ref{bsk}) with
\[
\sigma \sigma^*(x) = a(x,
u^{\varepsilon_n}(T^{\varepsilon_n}(\lambda_n), O_1)) = a(x,
{c}^1(\lambda) + \delta'))
\]
 will be in an arbitrarily small neighborhood of
$O_1$ at time $ T^{\varepsilon_n}(\lambda_n) -
T^{\varepsilon_n}(\mu_n)$ with probability which tends to one when
$\varepsilon_n \downarrow 0$.  By Lemma~\ref{mlml}, this remains
true if the constant
$u^{\varepsilon_n}(T^{\varepsilon_n}(\lambda_n), O_1)$ is replaced
by a function which is sufficiently close to this constant in
$D_1^\delta$, where $\delta$ is sufficiently small. Therefore, due
to (\ref{clll}) and Lemma~\ref{letwomk}, we can choose $\delta''$
sufficiently close to $\delta'$ so that $X^{
T^{\varepsilon_n}(\lambda_n),O_1,\varepsilon_n}_{T^{\varepsilon_n}(\lambda_n)
- T^{\varepsilon_n}(\mu_n)}$ will be in a small neighborhood of
$O_1$  with probability which tends to one when $\varepsilon_n
\downarrow 0$. With $\delta'$ and $\delta''$ thus fixed, we let
$\varepsilon_n \downarrow 0$ in (\ref{xpr}). The left hand side is
equal to ${c}^1(\lambda) + \delta'$, while the right hand side
tends to ${c}^1(\lambda) + \delta''$. This leads to a
contradiction which proves that (\ref{twhh}) holds, which in turn
implies that (\ref{toshow}) holds. The proof of (\ref{toshow2}) is
completely similar.

Note that the arguments used to prove (\ref{twhh}) also lead to
the following statement: for each $\lambda_0 > 0$
\begin{equation} \label{twhh3}
\limsup_{\varepsilon \downarrow 0} \sup_{\lambda' \in [\lambda_0,
\lambda]} u^\varepsilon(T^\varepsilon(\lambda'), O_1) \leq
\lim_{\lambda' \downarrow \lambda} {c}^1(\lambda'),
\end{equation}
now without assuming that $c^1$ is continuous at $\lambda$.
Similarly, for each $\lambda_0 > 0$
\begin{equation} \label{twhh4}
\liminf_{\varepsilon \downarrow 0} \inf_{\lambda' \in [\lambda_0,
\lambda]} u^\varepsilon(T^\varepsilon(\lambda'), O_2) \geq
\lim_{\lambda' \downarrow \lambda} {c}^2(\lambda').
\end{equation}
Let us show that if ${c}^1$ is continuous at $\lambda$, then
\begin{equation} \label{toshowx}
\liminf_{\varepsilon \downarrow 0} \inf_{x \in D_1^\delta}
u^\varepsilon(T^\varepsilon(\lambda), x) \geq \min({c}^1(\lambda),
\lim_{\lambda' \downarrow \lambda} {c}^2(\lambda')).
\end{equation}
Similarly, if  ${c}^2$ is continuous at $\lambda$, then
\begin{equation} \label{toshow2x}
\limsup_{\varepsilon \downarrow 0} \sup_{x \in D_2^\delta}
u^\varepsilon(T^\varepsilon(\lambda), x) \leq \max({c}^2(\lambda),
\lim_{\lambda' \downarrow \lambda} {c}^1(\lambda')).
\end{equation}
Due to Lemma~\ref{letwomk}, in order to prove (\ref{toshowx}), it
is sufficient to show that
\begin{equation} \label{toshowxz}
\liminf_{\varepsilon \downarrow 0}
u^\varepsilon(T^\varepsilon(\lambda), O_1) \geq
\min({c}^1(\lambda), \lim_{\lambda' \downarrow \lambda}
{c}^2(\lambda')).
\end{equation}
If (\ref{toshowxz}) fails, then for each $\lambda_0 > 0$ there is
$\delta'
> 0$ and a sequence $\varepsilon_n \downarrow 0$ such that
\begin{equation} \label{jji}
u^{\varepsilon_n}(T^{\varepsilon_n}(\lambda), O_1) <
{c}^1(\lambda) - \delta'.
\end{equation}
\begin{equation} \label{iio}
u^{\varepsilon_n}(T^{\varepsilon_n}(\lambda), O_1) <
\inf_{\lambda' \in [\lambda_0, \lambda]}
u^{\varepsilon_n}(T^{\varepsilon_n}(\lambda'), O_2) - \delta'.
\end{equation}
These two inequalities can not hold at the same time as follows
from Lemma 3.11 of \cite{FK}, where an analogue of (\ref{jji}) is
ruled out for the case of the initial-boundary value problem with
one equilibrium point inside the domain. Now the boundary
condition is replaced by the presence of the second equilibrium
point, but due to (\ref{iio}) the proof goes through without major
modifications. We have thus justified (\ref{toshowx}), and
(\ref{toshow2x}) is absolutely similar.


 Note that (\ref{toshow}), (\ref{toshow2}),
(\ref{toshowx}), and (\ref{toshow2x}) imply the statement of the
theorem for  $0< \lambda < \lambda_3$. Expressing the solution at
time $T^\varepsilon(\lambda)$ in terms of the solution at an
earlier time $T^\varepsilon(\lambda')$ (similarly to (\ref{xpr})),
we see that if
\[
\liminf_{\varepsilon \downarrow 0} \inf_{x \in D_1^\delta}
u^\varepsilon(T^\varepsilon(\lambda'), x) \leq
\limsup_{\varepsilon \downarrow 0} \sup_{x \in D_2^\delta}
u^\varepsilon(T^\varepsilon(\lambda'), x),
\]
then
\[
 \liminf_{\varepsilon \downarrow 0} \inf_{x \in D_1^\delta}
u^\varepsilon(T^\varepsilon(\lambda'), x) \leq
\liminf_{\varepsilon \downarrow 0} \inf_{x \in D_1^\delta \cup
D_2^\delta} u^\varepsilon(T^\varepsilon(\lambda), x) \leq
\]
\[
\leq \limsup_{\varepsilon \downarrow 0} \sup_{x \in D_1^\delta
\cup D_2^\delta} u^\varepsilon(T^\varepsilon(\lambda), x) \leq
\limsup_{\varepsilon \downarrow 0} \sup_{x \in D_2^\delta}
u^\varepsilon(T^\varepsilon(\lambda'), x).
\]
As follows from the definition of the functions
$\overline{c}^1(\lambda)$ and $\overline{c}^2(\lambda)$, this
allows us to extend the result to $\lambda \geq \lambda_3$.
 \qed
\\

\noindent {\bf Remark.} If $\lambda > \lambda_3$, then $
\overline{c}^1(\lambda) = \overline{c}^2(\lambda) =  c^*$.  It is
possible to show that the limit
\[
\lim_{\varepsilon \downarrow 0}
u^\varepsilon(T^{\varepsilon}(\lambda), x) = c^*
\]
is uniform in $(x, \lambda) \in B_{1/\delta} \times
[\overline{\lambda},\infty)$ for each $\overline{\lambda} >
\lambda_3$, where $ B_{1/\delta}$ is the ball of radius $1/\delta$
centered at the origin. Therefore, for each $\delta
> 0$ and $\overline{\lambda}
> \lambda_3$ there is $\varepsilon_0
> 0$ such that
\[
|u^\varepsilon(t, x) - c^*| \leq \delta
\]
whenever $\varepsilon \in (0,\varepsilon_0)$, $x \in B_{1/\delta}$
and $t \geq T^{\varepsilon}(\overline{\lambda})$.
\\

Let $X^{T^\varepsilon(\lambda),x,\varepsilon}_s$, $s \in
[0,T^\varepsilon(\lambda)]$, be the process defined in
(\ref{smtttnew})-(\ref{printe}). As follows from the large
deviation theory (see Chapter 6 of \cite{FW}), the distribution of
the random variable
$X^{T^\varepsilon(\lambda),x,\varepsilon}_{T^\varepsilon(\lambda)}$
will be concentrated near the points $O_1$ and $O_2$. From
Theorem~\ref{mtex1} and the representation (\ref{printe}) for the
solution, we obtain the following theorem.
\begin{theorem} \label{tenn}
Suppose that $g(O_1) < g(O_2)$. If the function
$\overline{c}^1(\lambda)$ is continuous at a point $\lambda \in
(0, \infty)$ and $x \in D_1$, then the distribution of the random
variable
$X^{T^\varepsilon(\lambda),x,\varepsilon}_{T^\varepsilon(\lambda)}$
converges to the measure $\mu_1^\lambda = a_1 \delta_{O_1} + a_2
\delta_{O_2}$, where the coefficients $a_1$ and $a_2$ can be found
from the equations $\overline{c}^1(\lambda) = a_1 g(O_1) + a_2
g(O_2)$, $a_1 + a_2 = 1$.

If the function $\overline{c}^2(\lambda)$ is continuous at a point
$\lambda \in (0, \infty)$ and $x \in D_2$, then the distribution
of the random variable
$X^{T^\varepsilon(\lambda),x,\varepsilon}_{T^\varepsilon(\lambda)}$
converges to the measure $\mu_2^\lambda = a_1 \delta_{O_1} + a_2
\delta_{O_2}$, where the coefficients $a_1$ and $a_2$ can be found
from the equations $\overline{c}^2(\lambda) = a_1 g(O_1) + a_2
g(O_2)$, $a_1 + a_2 = 1$.

If $\lambda \in (\lambda_3, \infty)$ and $x \in D$, then the
distribution of the random variable
$X^{T^\varepsilon(\lambda),x,\varepsilon}_{T^\varepsilon(\lambda)}$
converges to the measure $\mu^* = a_1 \delta_{O_1} + a_2
\delta_{O_2}$, where the coefficients $a_1$ and $a_2$ can be found
from the equations $c^* = a_1 g(O_1) + a_2 g(O_2)$, $a_1 + a_2 =
1$.
\end{theorem}

\section{Three equilibrium points without changes in the hierarchy
of cycles} \label{example2} In this section we  assume that there
are three asymptotically stable equilibrium points $O_1, O_2, O_3$
such that $g(O_1) \leq g(O_2) \leq g(O_3)$. For $c_1,c_2,c_3 \in
[g_{\rm min}, g_{\rm max}]$, let
\begin{equation} \label{tcs}
f_{c_1, c_2,c_3}(x) = c_1 \chi_{D_1}(x) + c_2 \chi_{D_2}(x) + c_3
\chi_{D_3}(x),~~x \in \mathbb{R}^d.
\end{equation}
Recall the definition of the hierarchy of cycles from
Section~\ref{defs}. We will assume that, for each
 choice of constants $c_i \in [g_{\rm
min},g_{\rm max}]$ in the function $\alpha = a(x,f_{c_1,
c_2,c_3}(x))$, Assumption~A holds and $O_1$ and $O_2$ form a cycle
$\Gamma$ of rank one. Consequently $O_1$, $O_2$ and $O_3$ form a
cycle of rank two for each
 choice of the constants.  Define
\[
M_{12}(c) =  {V}^{ a(\cdot, c)}_{O_1,O_2}, ~~~M_{21}(c) =  {V}^{
a(\cdot, c)}_{O_2,O_1},
\]
\[
M_{\Gamma3}(c) =  {V}^{ a(\cdot,
c)}_{\Gamma,O_3},~~~M_{3\Gamma}(c) = {V}^{ a(\cdot,
c)}_{O_3,\nu(O_3)}.
\]
Let
$\lambda_1 = M_{12}(g(O_1))$ and $\lambda_2 = M_{21}(g(O_2))$.
Define functions $c^1$ and $c^2$ by (\ref{hhnn1})
and~(\ref{hhnn2}), respectively. Let $\lambda_3 = \inf\{\lambda:
c^1(\lambda) \geq c^2(\lambda)\}$. Assume that at least one of the
functions $c^1$ and $c^2$ is continuous at $\lambda_3$. Let $c^* =
c^1(\lambda_3)$ if $c^1$ is continuous at $\lambda_3$ and $c^* =
c^2(\lambda_3)$ otherwise. Let $\overline{c}^1(\lambda) =
\min(c^1(\lambda), c^*)$ and $\overline{c}^2(\lambda) =
\max(c^2(\lambda), c^*)$, $\lambda < \lambda_3$. Let $\lambda_4 =
M_{\Gamma3}(c^*)$ and $\lambda_5 = M_{3\Gamma}(g(O_3))$.

\begin{figure}[htbp]
 \label{pic2b}
  \begin{center}
    \begin{psfrags}
     \includegraphics[height=6.5in, width= 3.5in,angle=90]{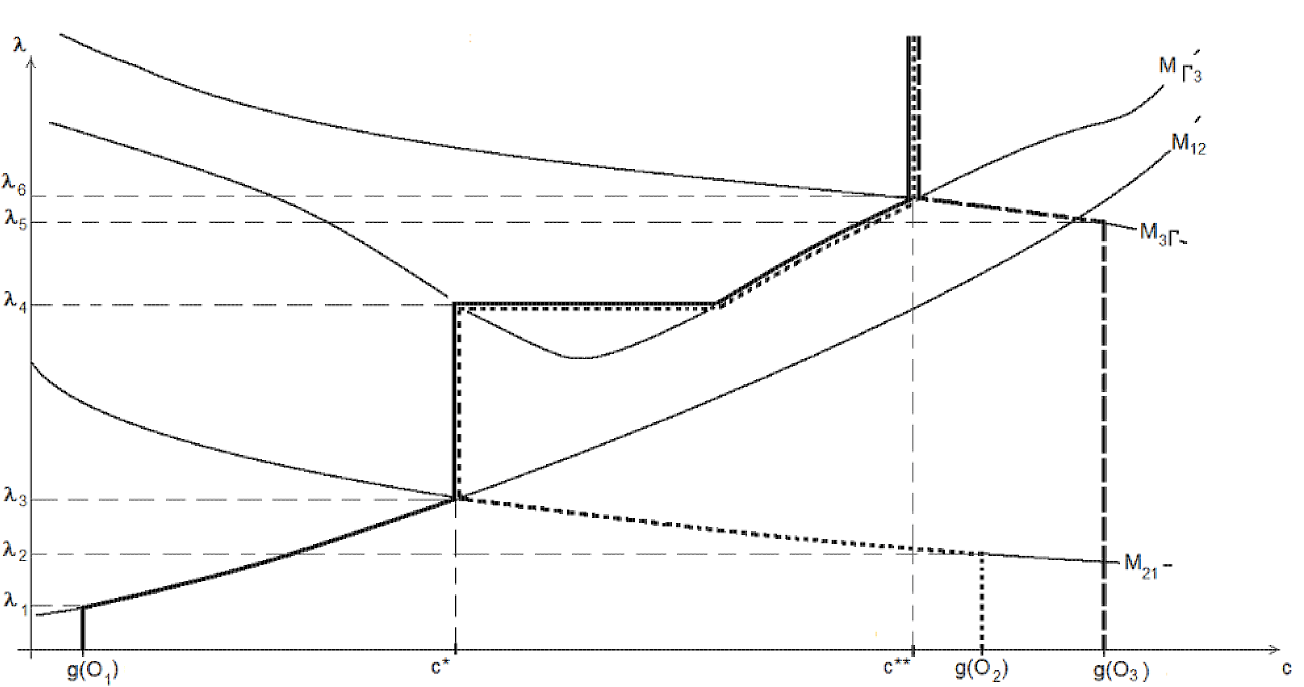}
    \end{psfrags}
    \caption{A case of three equilibrium points without changes in the hierarchy of cycles}
  \end{center}
\end{figure}

Let us assume that $\lambda_1 < \lambda_2 < \lambda_3 < \lambda_4
< \lambda_5 $ (see Figure 2). For $\lambda < \lambda_3$, the
behavior of the solution in $D_1$ and $D_2$ is still governed by
Theorem~\ref{mtex1}. For each $\lambda > \lambda_3$, the value of
$u^\varepsilon(T^\varepsilon(\lambda), x)$ will be nearly constant
on $D_1^\delta \cup D_2^\delta$, and we can treat the cycle
$\Gamma = \{O_1, O_2\}$ in the same way a single equilibrium was
treated in Section~\ref{example1}. Namely, let
\[
c^\Gamma(\lambda) =
\left\{\begin{array}{lll}c^*,~~~~~~~~~~~~~~~~~~~~~~~~~~~~~~~~~~~~~~~~~~~~~~~~~~~~~~~~~~~~~~~\,\lambda_3
\leq \lambda < \lambda_4,\\ \min \{g(O_3), \min\{c: c \in [c^*,
g(O_3)], M_{\Gamma3}(c) = \lambda \} \},~~\,\,\lambda \geq
\lambda_4,\end{array}\right.
\]
\[
c^3(\lambda) =
\left\{\begin{array}{lll}g(O_3),~~~~~~~~~~~~~~~~~~~~~~~~~~~~~~~~~~~~~~~~~~~~~~~~~~~~~~\,0
< \lambda < \lambda_5,\\ \max\{c^*, \max\{c: c \in [c^*, g(O_3)],
M_{3\Gamma}(c) = \lambda \} \},~~\,\,\lambda \geq
\lambda_5,\end{array}\right.
\]

%
%
%
%

 Define $\lambda_6 =
\inf\{\lambda > \lambda_3: c^\Gamma(\lambda) \geq c^3(\lambda)\}$.
Assume that $\lambda_5 < \lambda_6$ and that at least one of the
functions $c^\Gamma$ and $c^3$ is continuous at $\lambda_6$. Let
$c^{**} = c^\Gamma(\lambda_6)$ if $c^\Gamma$ is continuous at
$\lambda_6$ and $c^{**} = c^3(\lambda_6)$ otherwise. Define
$\overline{c}^1(\lambda) = \overline{c}^2(\lambda) =
\min(c^\Gamma(\lambda), c^{**})$,
$\lambda \geq \lambda_3$, and $\overline{c}^3(\lambda) =
\max(c^3(\lambda), c^{**})$, $\lambda
> 0$.

Having thus defined the functions $\overline{c}^i(\lambda)$, $i
=1,2,3$, for all $\lambda  > 0$, we can now state that for each
$\lambda > 0$ such that $\overline{c}^i$ is continuous at
$\lambda$ and every $\delta
> 0$, the limit
\[
\lim_{\varepsilon \downarrow 0}
u^\varepsilon(T^\varepsilon(\lambda), x) = \overline{c}^i(\lambda)
\]
is uniform in $x \in D_i^\delta$.

%
%
On Figure 2, the limits $ \lim_{\varepsilon \downarrow 0}
u^\varepsilon(T^\varepsilon(\lambda), x)$, as functions of
$\lambda$, for $x \in D_1^\delta$,  $ D_2^\delta$ and $D_3^\delta$
are depicted by thick, dotted and dashed lines, respectively.

\section{A general result for the case when the hierarchy of cycles does not change}
\label{general}
In this section we will suppose that, in addition to Assumption A,
the hierarchy of cycles and the equilibrium points $\nu(\Gamma)$
for each cycle $\Gamma$ of rank less than $R$ do not depend on the
choice of constants $c_i \in [g_{\rm min},g_{\rm max}]$ in the
function $\alpha = a(x,\sum_{i=1}^n c_i \chi_{D_i}(x))$.


We will say that a cycle $\Gamma$ is active for a given value of
$\lambda
> 0$ if $V^\alpha_{\Gamma, \nu(\Gamma)} < \lambda$. We will say that it
is engaged if $V^\alpha_{\Gamma, \nu(\Gamma)} = \lambda$ and
passive if $V^\alpha_{\Gamma, \nu(\Gamma)} > \lambda$. We will say
that a cycle $\Gamma_0$ is connected to a cycle $\Gamma$ by a
chain if there is a sequence of cycles $\Gamma_1,...,\Gamma_k$ and
equilibriums $O_1 \in \Gamma_1$,...,$O_k \in \Gamma_k$, $O_{k+1}
\in \Gamma$ such that $\Gamma_i$ are engaged or active  and
$O_{i+1} = \nu(\Gamma_i)$ for $0 \leq i \leq k$. The collection of
all the cycles that do not belong to $\Gamma$ and are connected to
$\Gamma$ by a chain will be called the cluster connected to
$\Gamma$. For each cycle $\Gamma$ of less than maximal rank  and
$c \in [g_{\rm min},g_{\rm max}]$, we define
\[
M_\Gamma(c) = V_{\Gamma, \nu(\Gamma)}^{a(\cdot, c)}
\]
and, for $\lambda > 0$ and  $c_2 \geq c_1$, define $
C(c_1,c_2,\lambda,\Gamma) = \min(c_2, \inf(c > c_1: M_\Gamma(c)
\geq \lambda))$. Similarly, if $c_2 \leq c_1$, define $
C(c_1,c_2,\lambda,\Gamma) = \max(c_2, \sup(c < c_1: M_\Gamma(c)
\geq \lambda))$.

In Figure 3 we have an example of a hierarchy of cycles with the
thick arrows between the actively connected cycles and the
corresponding equilibrium points. The dashed arrows are used for
the engaged cycles and the dotted arrows for the passively
connected cycles.

\begin{figure}[htbp]
 \label{pic2acc}
  \begin{center}
    \begin{psfrags}
     \includegraphics[height=5.8in, width= 3.5in,angle=90]{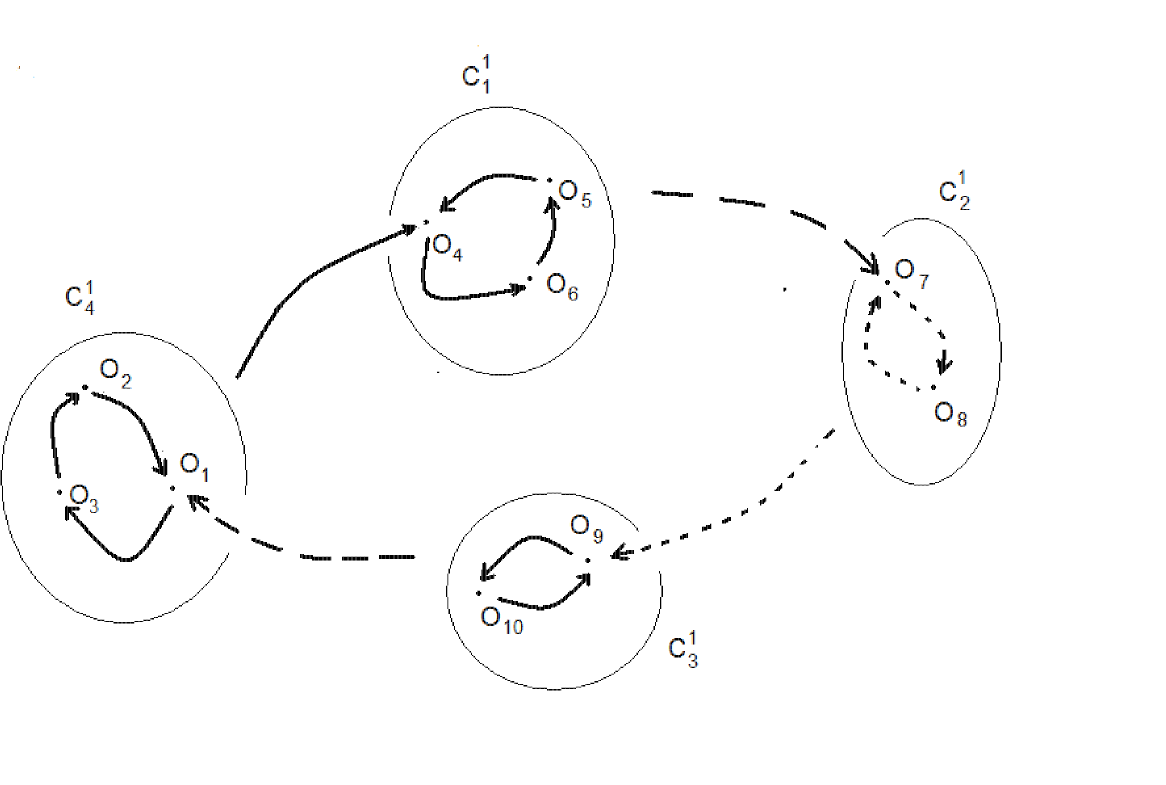}
    \end{psfrags}
    \caption{The hierarchy of cycles}
  \end{center}
\end{figure}

In order to describe the asymptotics of
$u^\varepsilon(T^\varepsilon(\lambda), x)$, we will define a
finite number of ``special" points $0 = \lambda_0 < \lambda_1 <
... < \lambda_{m} = \infty$.
 We claim that there
are functions $ \overline{c}^i(\lambda)$, $1 \leq i \leq n$, which
are continuous on each of the intervals $(\lambda_k,
\lambda_{k+1})$, $0 \leq k \leq m-1$, have one-sided limits as
$\lambda$ approaches the end points of the intervals, and are such
that the limits
\[
\lim_{\varepsilon \downarrow 0}
u^\varepsilon(T^\varepsilon(\lambda), x) = \overline{c}^i(\lambda)
\]
are uniform in $x \in D_i^\delta$ for each $\delta > 0$, $\lambda
\in \mathbb{R}^+ \setminus \Lambda$ with $\Lambda = \{\lambda_0,
\lambda_1,....,\lambda_m \}$. Moreover, neither of the cycles
changes its type (between passive, engaged and active) for
$\lambda \in (\lambda_k, \lambda_{k+1})$ and $\alpha(x) =
\lim_{\varepsilon \downarrow 0} a(x,
u^\varepsilon(T^\varepsilon(\lambda), x))$. We will use induction
on $k$ in order to define the  functions $\overline{c}^i(\lambda)$
and describe for each cycle whether it is passive, engaged or
active for $\lambda \in (\lambda_k, \lambda_{k+1})$ with
 $\alpha(x) = \lim_{\varepsilon \downarrow 0} a(x,
u^\varepsilon(T^\varepsilon(\lambda), x))$. In the process, we
will make several assumptions about the functions $M_\Gamma$.

%
%
%
%
%
%

Assuming that we have defined $\overline{c}^i(\lambda)$, let
\[
\lambda_\Gamma = \inf\{\lambda > 0: \overline{c}^i(\lambda')~~{\rm
does}~{\rm not}~{\rm depend}~ {\rm on}~ i~{\rm for}~\lambda' \geq
\lambda~{\rm and}~ O_i \in \Gamma \}.
\]
From  the inductive construction of the functions
$\overline{c}^i(\lambda)$ it will follow that $\lambda_\Gamma <
\infty$. Let $a_\Gamma = \lim_{\lambda \downarrow
\lambda_\Gamma}\overline{c}^i(\lambda)$, $O_i \in \Gamma$, and
$A_\Gamma = M_{\Gamma}(a_\Gamma)$. We assume that all $A_\Gamma$
are distinct and define
\[
\Lambda^1 = \{A_\Gamma,~{\rm rank}(\Gamma) < R \}.
\]

We assume that $M_\Gamma$ has a finite number of critical points
on $[g_{\rm min},g_{\rm max}]$ for each $\Gamma$ with ${\rm
rank}(\Gamma) < R$. Let $c^\Gamma_1,...,c^\Gamma_{k_\Gamma}$ be
all the local maxima of $M_\Gamma$. We assume that
$M_{\Gamma}(c^\Gamma_i)$ are distinct for all $\Gamma$  with ${\rm
rank}(\Gamma) < R$ and $i$. Define
\[
\Lambda^2 = \{M_{\Gamma}(c^\Gamma_i),~~{\rm rank}(\Gamma) < R,~1
\leq i \leq k_\Gamma \}.
\]

Let $\Gamma$ be a cycle of rank $r < R$, $\overline{\Gamma}$ the
cycle of rank $r +1$ that contains $\Gamma$, and $\Upsilon$ a
cycle that is contained in $\overline{\Gamma} \setminus \Gamma$.
Let $I_{\Gamma, \Upsilon} = \{c: M_\Gamma(c) = M_\Upsilon(c)\}$.
We assume that the sets $I_{\Gamma, \Upsilon}$ are finite and
$I_{\Gamma_1, \Upsilon_1} \cap  I_{\Gamma_2, \Upsilon_2} =
\emptyset$  unless $(\Gamma_1, \Upsilon_1) =
(\Gamma_2,\Upsilon_2)$ or $(\Gamma_1, \Upsilon_1) =
(\Upsilon_2,\Gamma_2)$. Define
\[
\Lambda^3 = \{M_\Gamma(c),~c \in   I_{\Gamma, \Upsilon},~~{\rm
rank}(\Gamma) < R, \Upsilon \subseteq \overline{\Gamma} \setminus
\Gamma \}.
\]
We assume that the numbers $M_\Gamma(a_{\Upsilon})$ are distinct
for all choices of cycles $\Gamma$ and $\Upsilon$ such that ${\rm
rank}(\Gamma) < R$, ${\rm rank}(\Upsilon) \leq {\rm
rank}({\Gamma})$ and $\nu(\Gamma) \in \Upsilon$. Define
\[
\Lambda^4 = \{M_\Gamma(a_{\Upsilon}),~~{\rm rank}(\Gamma) < R,
~{\rm rank}(\Upsilon) \leq {\rm rank}({\Gamma}),~ \nu(\Gamma) \in
\Upsilon \}.
\]

Finally, we assume that the sets $\Lambda^1$, $\Lambda^2$,
$\Lambda^3$ and $\Lambda^4$ do not intersect and define
\[
\Lambda = \{\lambda_0, \lambda_1,....,\lambda_m \} :=  \{0\} \cup
\Lambda^1 \cup \Lambda^2 \cup \Lambda^3 \cup \Lambda^4 \cup
\{\infty \},
\]
where we arrange $\lambda_k$ in the increasing order.

Below we will define $\overline{c}^i(\lambda)$ on the successive
intervals $(\lambda_k, \lambda_{k+1})$ using induction on $k$
while assuming that $\lambda_k$ are known. The above definition of
$\Lambda^1$ and $\Lambda^4$ in terms of $\overline{c}^i(\lambda)$
does not constitute a circular argument, since we could instead
define the pairs $(\lambda_{k+1}, \overline{c}^i(\lambda)~{\rm
for}~ \lambda \in (\lambda_k, \lambda_{k+1}))$ inductively.
Such an approach would lead to more complicated notations, though,
so we avoid it.

Let us proceed with the inductive definition of
$\overline{c}^i(\lambda)$.  For $\lambda \in (\lambda_0,
\lambda_1)$ all cycles are passive and $ \overline{c}^i(\lambda) =
g(O_i)$ for all~$i$.  Assuming that the types of the cycles and
the limits $q(O_i) = \lim_{\lambda \uparrow \lambda_k}
\overline{c}^i(\lambda)$ are known for $\lambda \in
(\lambda_{k-1}, \lambda_k)$ with some $0 <  k < m$, we will
describe the types of the cycles for $\lambda \in (\lambda_{k},
\lambda_{k+1})$ and specify the limits $s(O_i) = \lim_{\lambda
\downarrow \lambda_k} \overline{c}^i(\lambda)$. Then, assuming
that the types of the cycles are specified for $\lambda \in
(\lambda_{k}, \lambda_{k+1})$ and the values of $s(O_i)$ are
known, we will define the functions $\overline{c}^i(\lambda)$ for
$\lambda \in (\lambda_{k}, \lambda_{k+1})$.
 We distinguish a
number of cases depending on whether $\lambda_k$ belongs to
$\Lambda^1$, $\Lambda^2$, $\Lambda^3$ or $\Lambda^4$.

First, however, we describe the procedure for determining the
values of $s(O_i)$ for $O_i$ which belong to a cluster.

{\bf Determining the values of $s(O_i)$ and the types of cycles
within a cluster.} Suppose that we have defined $s(O_i) =
\lim_{\lambda \downarrow \lambda_k} \overline{c}^i(\lambda)$ for
all $O_i$ that belong to a cycle $\Gamma$. Consider the cluster of
cycles that are connected to $\Gamma$ for $\lambda \in
(\lambda_{k-1}, \lambda_{k})$. For each cycle $\Gamma'$ in the
cluster, we will define the values of $s(O)$ for $O \in \Gamma'$
and specify its type for $\lambda \in (\lambda_{k},
\lambda_{k+1})$.

First assume that $\nu(\Gamma')=O_i \in \Gamma$ for $\lambda \in
(\lambda_{k-1}, \lambda_{k})$. It will follow from the inductive
construction that $q({O'}) = q(O'')$ if $O', O'' \in \Gamma'$. Let
$q(\Gamma') = q({O'})$. For $O \in \Gamma'$, we define $s(O) =
C(q(\Gamma'),s(O_i),\lambda_k,\Gamma')$.

For any cycle $\Gamma''$  such that $\nu(\Gamma'') \in \Gamma'$,
we can similarly  determine the values of $s(O)$ for $O \in
\Gamma''$. Continuing this procedure inductively, we define the
values of $s(O)$ when $O$ belongs to either of the cycles from the
cluster. A cycle $\Gamma'$ from the cluster will be engaged for
$\lambda \in (\lambda_{k}, \lambda_{k+1})$ if $ \lambda_k =
M_{\Gamma'}(s(O))$ for $O \in \Gamma'$ and active if $ \lambda_k
> M_{\Gamma'}(s(O))$ for $O \in \Gamma'$.

 {\bf Case 1.} Assume that $\lambda_k \in \Lambda_1$. Let $\Gamma$ be
 such that $\lambda_k = A_\Gamma$. For $O_i \in \Gamma$, we define
$s(O_i) = C(q(O_i),q(\nu(\Gamma)),\lambda_k,\Gamma)$. The cycle
 will be engaged for
$\lambda \in (\lambda_{k}, \lambda_{k+1})$ if $ \lambda_k =
M_{\Gamma}(s(O))$ for $O \in \Gamma$ and active if $ \lambda_k >
M_{\Gamma}(s(O))$ for $O \in \Gamma$.

The types of cycles  that belong the the cluster connected to
$\Gamma$ for $\lambda \in (\lambda_{k-1}, \lambda_k)$, and the
values of $s(O_j)$ for the equilibrium points in those cycles are
determined according to the procedure described above. For the
remaining equilibrium points $O$, we define $s(O) = q(O)$. The
remaining cycles don't change type.

{\bf Case 2.} Assume that $\lambda_k \in \Lambda_2$. Let $c$ be
the local maximum of a cycle $\Gamma$ such that $M_\Gamma(c) =
\lambda_k$. If $\Gamma$ was not engaged for $\lambda \in
(\lambda_{k-1}, \lambda_k)$ or if $q(O) \neq c$ for some $O \in
\Gamma$, then we define $s(O) = q(O)$ for all the equilibrium
points, and all the cycles have the same type on $(\lambda_{k},
\lambda_{k+1})$ as on $(\lambda_{k-1}, \lambda_k)$.

If $\Gamma$ was engaged and $q(O) = c$ for $O \in \Gamma$, then
for $O_i \in \Gamma$, we define $s(O_i) =
C(q(O_i),q(\nu(\Gamma)),\lambda_k,\Gamma)$. The cycle
 will be engaged for
$\lambda \in (\lambda_{k}, \lambda_{k+1})$ if $ \lambda_k =
M_{\Gamma}(s(O))$ for $O \in \Gamma$ and active if $ \lambda_k >
M_{\Gamma}(s(O))$ for $O \in \Gamma$.

The types of cycles  that belong the the cluster connected to
$\Gamma$ for $\lambda \in (\lambda_{k-1}, \lambda_k)$, and the
values of $s(O_j)$ for the equilibrium points in those cycles are
determined according to the procedure described above. For the
remaining equilibrium points $O$, we define $s(O) = q(O)$. The
remaining cycles don't change type.

{\bf Case 3.} Assume that $\lambda_k \in \Lambda_3$. Let $\Gamma$
be a cycle of rank $r < R$, $\overline{\Gamma}$ the cycle of rank
$r +1$ that contains $\Gamma$, and $\Upsilon$ a cycle that is
contained in $\overline{\Gamma} \setminus \Gamma$. Suppose that
$c$ is such that $ M_\Gamma(c) = M_\Upsilon(c)$ and $\lambda_k =
M_\Gamma(c)$.

We define $s(O) = q(O)$ for all the equilibrium points. All the
cycles, other than perhaps $\Gamma$ and $\Upsilon$, have the same
type on $(\lambda_{k}, \lambda_{k+1})$ as on $(\lambda_{k-1},
\lambda_k)$. To determine the type of cycles $\Gamma$ and
$\Upsilon$ on $(\lambda_{k}, \lambda_{k+1})$, we examine several
cases.

(a) If $q(O) = c$ for all $O \in \Gamma \cup \Upsilon$, $\Gamma$
and $\Upsilon$ were engaged, $\Upsilon$ was connected to $\Gamma$
by a chain that contained only active cycles (other than
$\Upsilon$ itself) and $\Gamma$ was connected to $\Upsilon$ by a
chain that contained only active cycles (other than $\Gamma$
itself),  then $\Gamma$ and $\Upsilon$ becomes active.

(b) If $q(O) = c$ for all $O \in \Gamma \cup \Upsilon$, $\Gamma$
was connected to $\Upsilon$ by a chain that contained only active
cycles (other than $\Gamma$ itself), but $\Upsilon$ was not
connected to $\Gamma$ by a chain that contained only active cycles
(other than $\Upsilon$ itself), and $\Upsilon$ was not passive,
then $\Gamma$ becomes active on $(\lambda_{k}, \lambda_{k+1})$ if
it was engaged on $(\lambda_{k-1}, \lambda_{k})$ and becomes
engaged if it was active. The type of $\Upsilon$ stays the same.

(c) the same as (b) with $\Gamma$ and $\Upsilon$ interchanged.

(d) If none of the cases (a)-(c) applies, then $\Gamma$ and
$\Upsilon$ have the same types on $(\lambda_{k}, \lambda_{k+1})$
as on $(\lambda_{k-1}, \lambda_k)$.

{\bf Case 4.} Assume that $\lambda_k \in \Lambda_4$. Let
$\lambda_k = M_\Gamma(a_{\Upsilon})$, where  cycles $\Gamma$ and
$\Upsilon$ are such that ${\rm rank}(\Gamma) < R$, ${\rm
rank}(\Upsilon) \leq {\rm rank}({\Gamma})$ and $\nu(\Gamma) \in
\Upsilon$. We define $s(O) = q(O)$ for all the equilibrium points.
All the cycles, other than perhaps $\Gamma$, have the same type on
$(\lambda_{k}, \lambda_{k+1})$ as on $(\lambda_{k-1}, \lambda_k)$.

The cycle $\Gamma$ becomes active if it was engaged on
$(\lambda_{k-1}, \lambda_k)$, $q(O) = a_\Upsilon$ for all $O \in
\Gamma$ and $M_\Gamma(a_\Upsilon) < A_\Upsilon$. Otherwise,
$\Gamma$  has the same type on $(\lambda_{k}, \lambda_{k+1})$ as
on $(\lambda_{k-1}, \lambda_k)$.
\\

Now let us define the functions $\overline{c}^i(\lambda)$ on
$(\lambda_k, \lambda_{k+1})$ assuming that the values of $s(O_i)$
and the cycle types are known. For an equilibrium point $O_i$, we
identify the cycle $\Gamma$ with the smallest possible rank $r$
such that $O_i \in \Gamma$ and the values of $s(O_j)$, $O_j \in
\Gamma$, are not all the same. If no such cycle exists, that is if
$s(O_j)$, $ 1 \leq j \leq n$, does not depend on $j$, then we
define $\overline{c}^i(\lambda) = s(O_i)$ for $\lambda >
\lambda_k$.

Assuming that such a cycle $\Gamma$ exists, let
$\Gamma_1,...,\Gamma_l$ be the cycles of rank $r-1$ which comprise
$\Gamma$, and let $O \in \Gamma_1$. Here we number the cycles in
such a way that $\mathcal{N}(\Gamma_1) =
\Gamma_2$,...,$\mathcal{N}(\Gamma_l) = \Gamma_1$. Take the least
$j$ such that $\Gamma_j$ is either passive or engaged (it can not
happen that all the cycles $\Gamma_1,...,\Gamma_l$ are active,
since then all the values of $s(O)$, $O \in \Gamma$, would be the
same, as follows from the inductive construction above). If
$\Gamma_j$ is passive, we define $\overline{c}^i(\lambda) =
s(O_i)$ for $\lambda \in (\lambda_k, \lambda_{k+1})$. If
$\Gamma_j$ is engaged, we define $\overline{c}^i(\lambda) =
C(r(O_i),\zeta,\lambda,\Gamma_j)$ for $\lambda \in (\lambda_k,
\lambda_{k+1})$, where $\zeta = +\infty$ if $M_{\Gamma_j}$ is
locally increasing at $r(O_i)$ and $\zeta = -\infty$ if
$M_{\Gamma_j}$ is locally decreasing at $r(O_i)$.

We can now summarize the above discussion.
\begin{theorem}
 Suppose that Assumption A holds and
the hierarchy of cycles and the equilibrium points $\nu(\Gamma)$
for each cycle $\Gamma$ of rank less than $R$ do not depend on the
choice of the constants $c_i \in [g_{\rm min},g_{\rm max}]$ in the
function $\alpha = a(x,\sum_{i=1}^n c_i \chi_{D_i}(x))$. Also
suppose that the above assumptions on the sets $\Lambda^1$,
$\Lambda^2$, $\Lambda^3$ and $\Lambda^4$ hold.

Then the limits
\[
\lim_{\varepsilon \downarrow 0}
u^\varepsilon(T^\varepsilon(\lambda), x) = \overline{c}^i(\lambda)
\]
are uniform in $x \in D_i^\delta$ for each $\delta > 0$, $\lambda
\in \mathbb{R}^+ \setminus \Lambda$, where the functions
$\overline{c}^i(\lambda)$ were defined via the inductive procedure
above.
\end{theorem}

\section{Example of a change in the hierarchy of cycles}
\label{change} As in Section~\ref{example2}, we assume that there
are three equilibrium points $O_1, O_2, O_3$.
For each $c_1,c_2,c_3 \in [g_{\rm min},
g_{\rm max}]$, the function $ f_{c_1, c_2,c_3}$ is defined
by~(\ref{tcs}). We will assume that the hierarchy of cycles for
$\alpha = a(x,f_{c_1, c_2,c_3}(x))$ depends only on $c_2$. This is
the case, for example, if $d = 1$ and $O_1 < O_2 < O_3$.
 More
precisely, suppose that there is $\overline{c} \in (g_{\rm min},
g_{\rm max})$ such that Assumption A holds for each choice of the
constants  $c_i \in [g_{\rm min}, g_{\rm max}]$ such that $c_2
\neq \overline{c}$. We assume that $O_1$ and $O_2$ form a cycle
$\Gamma' = \{O_1, O_2 \}$ of rank one when $c_2 < \overline{c}$,
while $O_2$ and $O_3$ form a cycle $\Gamma'' = \{O_2, O_3 \}$ of
rank one when $c_2
> \overline{c}$.

As before, we will identify a number of ``special" points
$\lambda_k$ and describe the asymptotic behavior of
$u^\varepsilon(T^\varepsilon(\lambda), x)$ for $\lambda \in
(\lambda_{k}, \lambda_{k+1})$ and $x \in D^\delta_i$, $i =1,2,3$.
In the process, we will make various assumptions about the
quasi-potential that will be specific to the example at hand.

In our example we assume that $g(O_1) \leq g(O_2) \leq
\overline{c} \leq g(O_3)$. Define
\[
M_{12}(c) =  {V}^{ a(\cdot, c)}_{O_1,O_2}, ~M_{21}(c) =  {V}^{
a(\cdot, c)}_{O_2,O_1}, ~M_{\Gamma' 3}(c) =  {V}^{ a(\cdot,
c)}_{\Gamma',O_3},~~c \in [g_{\rm min}, \overline{c});
\]
\[
M_{32}(c) =  {V}^{ a(\cdot, c)}_{O_3,O_2}, ~M_{23}(c) =  {V}^{
a(\cdot, c)}_{O_2,O_3}, ~M_{\Gamma'' 1}(c) =  {V}^{ a(\cdot,
c)}_{\Gamma'',O_1},~ M_{1 \Gamma''} (c) = {V}^{ a(\cdot, c)}_{O_1,
\nu(O_1)},~~c \in ( \overline{c}, g_{\rm max}].
\]

Let  $\lambda_1 = M_{12}(g(O_1))$ and $\lambda_2 =
M_{21}(g(O_2))$. Define functions $c^1$ and $c^2$ by (\ref{hhnn1})
and (\ref{hhnn2}), respectively. Let $\lambda_3 = \inf\{\lambda:
c^1(\lambda) \geq c^2(\lambda)\}$. Assume that at least one of the
functions $c^1$ and $c^2$ is continuous at $\lambda_3$. Let $c^* =
c^1(\lambda_3)$ if $c^1$ is continuous at $\lambda_3$ and $c^* =
c^2(\lambda_3)$ otherwise. Let $\overline{c}^1(\lambda) =
\min(c^1(\lambda), c^*)$ and $\overline{c}^2(\lambda) =
\max(c^2(\lambda), c^*)$, $\lambda < \lambda_3$. Let $\lambda_4 =
M_{\Gamma' 3}(c^*)$, $\lambda_5 = \sup_{c \in [c^*, \overline{c})}
M_{\Gamma' 3}({c})$ and $\lambda_6 = M_{32}(g(O_3))$.

\begin{figure}[htbp]
 \label{pic2c}
  \begin{center}
    \begin{psfrags}
     \includegraphics[height=6.5in, width= 3.5in,angle=90]{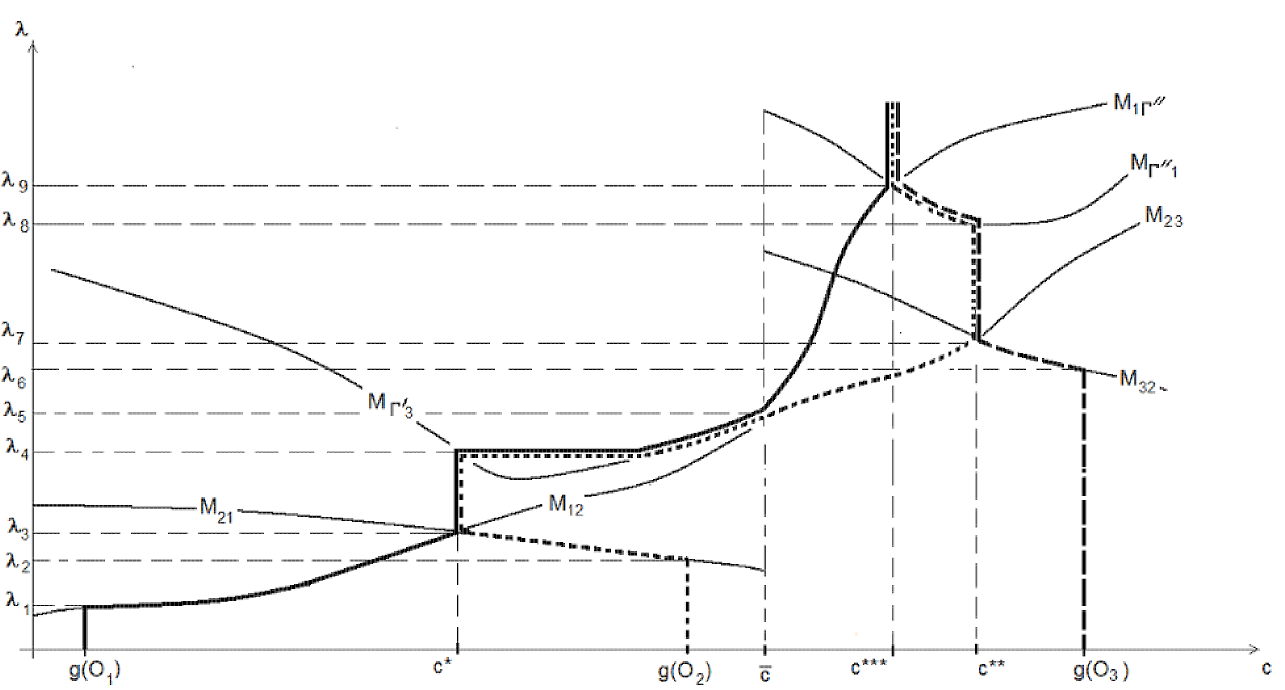}
    \end{psfrags}
    \caption{A case of three equilibrium when the hierarchy of cycles changes}
  \end{center}
\end{figure}

Let us assume that $\lambda_1 < \lambda_2 < \lambda_3 < \lambda_4
< \lambda_5 < \lambda_6$ (see Figure 4).
Define
\[
\overline{c}^1(\lambda) = \overline{c}^2(\lambda) =
c^{\Gamma'}(\lambda) =
\left\{\begin{array}{lll}c^*,~~~~~~~~~~~~~~~~~~~~~~~~~~~~~~~~~~~~~~~~~~~~~~~~~~~~~\,\lambda_3
\leq \lambda < \lambda_4,\\  \min\{c: c \in [c^*, \overline{c}),
M_{\Gamma'3}(c) = \lambda \},~~~~~~~~~~~~~~\,\,\lambda_4 \leq
\lambda < \lambda_5,\end{array}\right.
\]
In order to
formulate the results on the asymptotics of
$u^\varepsilon(T^\varepsilon(\lambda), x)$ for $\lambda >
\lambda_5$, we need the functions $d^2(\lambda)$ and
$c^3(\lambda)$ defined as follows:
\[
d^2(\lambda) = \min \{g(O_3), \min\{c: c \in [\overline{c},
g(O_3)], M_{23}(c) = \lambda \} \},~~\lambda \geq \lambda_5,
\]
\[
c^3(\lambda) =
\left\{\begin{array}{lll}g(O_3),~~~~~~~~~~~~~~~~~~~~~~~~~~~~~~~~~~~~~~~~~~~~~~~~~~~~~~\,0
< \lambda < \lambda_6,\\ \max\{\overline{c}, \max\{c: c \in
[\overline{c}, g(O_3)], M_{32}(c) = \lambda \} \},~~~~\,\,\lambda
\geq \lambda_6.\end{array}\right.
\]
 Let $\lambda_7 = \inf\{\lambda:
d^2(\lambda) \geq c^3(\lambda)\}$. Assume that $\lambda_6 <
\lambda_7$ and at least one of the functions $d^2$ and $c^3$ is
continuous at $\lambda_7$. Let $c^{**} = d^2(\lambda_7)$ if $d^2$
is continuous at $\lambda_7$ and $c^{**} = c^3(\lambda_7)$
otherwise. Let $\lambda_8 = M_{\Gamma'' 1} (c^{**})$ and assume
that $\lambda_7 <  \lambda_8$. Define $\overline{c}^2(\lambda) =
\min(d^2(\lambda), c^{**})$, $\lambda_5 \leq \lambda < \lambda_8$,
and $\overline{c}^3(\lambda) = \max(c^3(\lambda), c^{**})$, $0 <
\lambda  < \lambda_8$.

 Let
\[ d^1(\lambda) = \min
\{c^{**}, \min\{c: c \in [\overline{c}, c^{**}], M_{1 \Gamma''}(c)
= \lambda \} \},~~\lambda \geq \lambda_5,
\]
\[
c^{\Gamma''}(\lambda) = \max\{\overline{c}, \max\{c: c \in
[\overline{c}, c^{**}], M_{\Gamma'' 1}(c) = \lambda \}
\},~~\lambda \geq \lambda_8.
\]
Let $\lambda_9 = \inf\{\lambda: d^1(\lambda) \geq
c^{\Gamma''}(\lambda)\}$. Assume that $\lambda_8 < \lambda_9$ and
at least one of the functions $d^1$ and $c^{\Gamma''}$ is
continuous at $\lambda_9$. Let $c^{***} = d^1(\lambda_9)$ if $d^1$
is continuous at $\lambda_9$ and $c^{***} =
c^{\Gamma''}(\lambda_9)$ otherwise.  Define
$\overline{c}^1(\lambda) = \min(d^1(\lambda), c^{***})$,
$\lambda_5 \leq \lambda$ and $\overline{c}^2(\lambda) =
\overline{c}^3(\lambda) = \max(c^{\Gamma''}(\lambda), c^{***})$,
$\lambda_8 \leq \lambda$.

Having thus defined the functions $\overline{c}^i(\lambda)$, $i
=1,2,3$, for all $\lambda  > 0$, we can now state that for each
$\lambda > 0$ such that $\overline{c}^i$ is continuous at
$\lambda$ and every $\delta
> 0$, the limit
\[
\lim_{\varepsilon \downarrow 0}
u^\varepsilon(T^\varepsilon(\lambda), x) = \overline{c}^i(\lambda)
\]
is uniform in $x \in D_i^\delta$.

On Figure 4, the limits $ \lim_{\varepsilon \downarrow 0}
u^\varepsilon(T^\varepsilon(\lambda), x)$, as functions of
$\lambda$, for $x \in D_1^\delta$,  $ D_2^\delta$ and $D_3^\delta$
are depicted by thick, dotted and dashed lines, respectively.
\\
\\
\\

\noindent {\bf \large Acknowledgements}: While working on this
article, M. Freidlin was supported by NSF grants DMS-0803287 and
DMS-0854982 and L. Koralov was supported by NSF grants DMS-0706974
and DMS-0854982.
\\
\\

\end{document}